\documentclass[10pt]{amsart}
\usepackage{amssymb}\usepackage{euscript}
\usepackage[all]{xy}

\usepackage{amsmath}

\newtheorem{theorem}{Theorem}
\newtheorem{lemma}[theorem]{Lemma}
\newtheorem{proposition}[theorem]{Proposition}
\newtheorem{corollary}[theorem]{Corollary} 
\newtheorem{definition}[theorem]{Definition}

\theoremstyle{definition}  
\newtheorem{example}[theorem]{Example}

\newtheorem{remark}[theorem]{Remark}

\newcommand{\nc}{\newcommand}
\newcommand{\on}{\operatorname}   

\newcommand{\ul}{\underline}

\newcommand{\kk}{{\bf k}}
\newcommand{\ZZ}{{\mathbb Z}}

\newcommand{\QQ}{{\mathbb Q}}

\renewcommand{\t}{{\mathfrak{t}}}
\newcommand{\g}{{\mathfrak{g}}}

\nc{\wh}{\widehat}\nc{\wt}{\widetilde}

\newcommand{\ben}{\begin{enumerate}}
\newcommand{\een}{\end{enumerate}}

\newcommand{\cO}{{\mathcal O}}

\newcommand{\cU}{{\mathcal U}}
\newcommand{\cC}{{\mathcal C}}

\newcommand{\cX}{{\mathcal X}}

\newcommand{\cD}{{\mathcal D}}

\newcommand{\cG}{{\mathcal G}}

\nc{\G}{{\mathfrak{G}}}
\nc{\HH}{{\mathfrak H}}

\newcommand{\p}{{\mathfrak{p}}}

\newcommand{\f}{{\mathfrak{f}}}

\begin{document}

\title[Half-balanced categories and Teichm\"uller groupoids in 
genus zero]{Half-balanced braided monoidal categories and Teichm\"uller 
groupoids in genus zero}

\author{Benjamin Enriquez}
\address{Institut de Recherche Math\'ematique Avanc\'ee, UMR 7501, 
Universit\'e de Strasbourg et CNRS, 
7 rue Ren\'e Descartes, 67000 Strasbourg, France}
\email{enriquez@math.unistra.fr}

\maketitle

\begin{abstract} 
We introduce the notions of a half-balanced braided monoidal category and of
its contraction. These notions give rise to an explicit description of 
the action of the Galois group of $\QQ$ on Teichm\"uller groupoids in 
genus 0, equivalent to that of L.~Schneps. We also show that a prounipotent 
version of this action is equivalent to a graded action.
\end{abstract}

\section*{Introduction and main results}

Let $M_{g,n}^\QQ$ be the moduli space of curves of genus $g$ with $n$
marked points. Its fundamental groupoid with respect to the set of 
maximally degenerate curves is called the Teichm\"uller groupoid $T_{g,n}$.
One of the main features of Grothendieck's geometric approach to the 
Galois group $G_\QQ$ of $\QQ$ is the study of its action on the profinite completions 
$\widehat T_{g,n}$; according to this philosophy, $G_\QQ$ could be 
characterized as the group of automorphisms of the tower of all 
$\widehat T_{g,n}$, compatible with natural operations, such as the 
Knudsen morphisms. It is therefore important to describe explicitly the 
action of $G_\QQ$ on the collection of all the $\widehat T_{0,n}$. 
Such a description was obtained in \cite{Sch}. More precisely, an explicit 
profinite group $\wh{\on{GT}}$ was introduced in \cite{Dr}, together with a 
morphism $G_\QQ\to\widehat{\on{GT}}$. The following was then proved in
\cite{Sch}:

\begin{theorem} \label{thm:1}
There exists a morphism $\widehat{\on{GT}}\to\on{Aut}(\widehat T_{0,n})$, 
such that the morphism $G_\QQ\to\on{Aut}(\widehat T_{0,n})$ factors
as $G_\QQ\to\widehat{\on{GT}}\to\on{Aut}(\widehat T_{0,n})$. 
\end{theorem}

The first purpose of this paper is to present a variant of the proof of
\cite{Sch}. This variant relies on the notion of a half-balanced braided 
monoidal category (b.m.c.), which appeared implicitly recently in \cite{ST} 
and is here made explicit. We introduce the notion of a (half-)balanced 
contraction of such a category $\cC$: it consists of a functor $\cC\to\cO$, 
satisfying certain properties. Whereas a balanced b.m.c.~gives rise to
representations of the framed braid group on the plane $\tilde B_n$ 
(for $n\geq 0$), which is an abelian extension of the braid group $B_n$, 
a (half-)balanced contraction gives rise to representations of quotients 
of $\tilde B_n$. This quotient is an abelian extension of the quotient 
$B_n/Z(B_n)$ of $B_n$ by its center in the case of a balanced contraction, and 
is an abelian extension of the mapping class group in genus zero $\Gamma_{0,n}$ 
(another  quotient of $B_n$) in the case of a half-balanced contraction. 

To each set $S$, 
we associate an object $\widehat{\bf PaB}_S^{hbal}\to\widehat{\bf PaDih}_S$ 
in the category whose objects are contractions of profinite half-balanced 
b.m.~categories, enjoying universal properties. These contractions may be viewed 
as the analogues of the universal b.m.~categories appearing in \cite{JS}. We show that 
the action of $\widehat{\on{GT}}$ on such categories may be lifted to the half-balanced 
setup. This defines in particular an action of $\widehat{\on{GT}}$ on 
$\widehat{\bf PaDih}_S$, from which it is is easy to derive an action of 
$\widehat{T}_{0,n}$.  

The above profinite theory admits a prounipotent version. The group 
$\widehat{\on{GT}}$ and the Teichm\"uller groupoid $\widehat{T}_{0,n}$
admit proalgebraic versions $\kk\mapsto \on{GT}(\kk),T_{0,n}(\kk)$, 
where $\kk$ is a $\QQ$-ring. We then have morphisms $\widehat{\on{GT}}
\to\on{GT}(\QQ_l)$, $\widehat T_{0,n}\to(T_{0,n})_l\to T_{0,n}(\QQ_l)$, 
where $l$ is a prime number and $(T_{0,n})_l$ is the pro-$l$ completion of 
$T_{0,n}$. We construct a group scheme $\ul{\on{Aut}T_{0,n}}(-)$,
together with a morphism $\ul{\on{Aut}T_{0,n}}(\kk)\to
\on{Aut}(T_{0,n}(\kk))$, a group scheme morphism 
$\on{GT}(-)\to\ul{\on{Aut}T_{0,n}}(-)$, and a group 
$\on{Aut}((T_{0,n})_l,T_{0,n}(\QQ_l))$, equipped with morphisms 
$$
\on{Aut}((T_{0,n})_l) \leftarrow 
\on{Aut}((T_{0,n})_l,T_{0,n}(\QQ_l)) \to 
\ul{\on{Aut}T_{0,n}}(\QQ_l). 
$$

\begin{theorem} \label{thm:2}
The morphism $G_{\QQ}\to \on{Aut}((T_{0,n})_{l})$ factors as 
$G_{\QQ}\to \on{Aut}((T_{0,n})_l,T_{0,n}(\QQ_l))\to \on{Aut}((T_{0,n})_{l})$, and 
there exists a morphism $\on{GT}(-)\to \ul{\on{Aut}T_{0,n}}(-)$, such that the following 
diagram commutes 
$$
\xymatrix{
G_{\QQ} \ar[r]\ar[dr]\ar[dd] & \on{Aut}((T_{0,n})_{l}) \\
 & \on{Aut}((T_{0,n})_l,T_{0,n}(\QQ_l))
 \ar[u]\ar[d]\\
 \on{GT}(\QQ_{l})\ar[r] & \ul{\on{Aut}T_{0,n}}(\QQ_{l})}$$
\end{theorem}

We say that an algebraic (resp., prounipotent) group over $\QQ$ is graded iff 
its Lie algebra is graded by $\ZZ_{\geq 0}$ (resp., by $\ZZ_{>0}$). 
We say that a groupoid $\cG$ is graded prounipotent if for any $s\in \on{Ob}\cG$, 
$\on{Aut}_{\cG}(s)$ is graded prounipotent. In \cite{Dr}, a graded 
$\QQ$-algebraic group $\on{GRT}(-)$ was constructed, together with an 
isomorphism $\on{GT}(-)\to\on{GRT}(-)$.  

\begin{theorem} \label{thm:3}
There exists a graded prounipotent groupoid $T_{0,n}^{gr}(-)$ and a graded 
morphism $\on{GRT}(-)\to \ul{\on{Aut}T_{0,n}^{gr}}(-)$, such that the 
diagram $\begin{matrix}\on{GT}(-)&\to&\ul{\on{Aut}T_{0,n}}(-) \\
\downarrow &&\downarrow \\
\on{GRT}(-)&\to&\ul{\on{Aut}T_{0,n}^{gr}}\end{matrix}$ commutes. 
\end{theorem}

\section{Teichm\"uller groupoids in genus 0} \label{sec:1}

\subsection{Quotient categories}

Let $\cC$ be a small category and let $G$ be a group. We define an action 
of $G$ on $\cC$ as the data of: (a) a group morphism $G\to \on{Perm}(\on{Ob}\cC)$, 
(b) for any $g\in G$, an assignment $\on{Ob}\cC\in X\mapsto i_X^{g}\in 
\on{Iso}_{\cC}(X,gX)$, such that $i_{X}^{gh} = i_{gX}^{h}i_{X}^{g}$. 

We then get a group morphism $G\to \on{Aut}\cC = \{$autofunctors of $\cC\}$, 
where the autofunctor induced by $g\in G$ is the action of $g$ at the level of objects,
and $g\phi := i_{Y}^{g}\phi(i_{X}^{g})^{-1}$ for $\phi\in\on{Hom}_{\cC}(X,Y)$.

\begin{lemma}
1) For any $\alpha,\beta\in(\on{Ob}\cC)/G$, there is a unique action of $G\times G$ 
on $\cX(\alpha,\beta):=
\sqcup_{X\in\alpha,Y\in\beta}\on{Hom}_{\cC}(X,Y)$, such that
$(g,h)\on{Hom}_{\cC}(X,Y) = \on{Hom}_{\cC}(gX,hY)$ and 
$(g,h)\phi = i_{Y}^{h}\phi(i_{X}^{g})^{-1}$. 

2) Set $\cX(X,\beta):= \sqcup_{Y\in\beta}\on{Hom}_{\cC}(X,Y)$, 
$\cX(\alpha,Y):= \sqcup_{X\in\alpha}\on{Hom}_{\cC}(X,Y)$, then 
$G$ acts on these sets (by permutation of $\beta$ in the first case and of $\alpha$ 
in the second 
one) and we have a well-defined map $\cX(X,\beta)^{G}\times \cX(\beta,Z)^{G}
\to \on{Hom}_{\cC}(X,Z)$ compatible all the maps $\on{Hom}_{\cC}(X,Y)
\times \on{Hom}_{\cC}(Y,Z)\to \on{Hom}_{\cC}(X,Z)$. Taking the product of 
these maps over $X\in\alpha$, $Z\in\gamma$ and further the quotient
by $G\times G$, we obtain a map $\cX(\alpha,\beta)^{G\times G}
\times\cX(\beta,\gamma)^{G\times G}
\to \cX(\alpha,\gamma)^{G\times G}$, which is associative. 
\end{lemma}

The proof is straightforward.  We then define the quotient category $\cC/G$ by 
$\on{Ob}(\cC/G):= (\on{Ob}\cC)/G$ and $(\cC/G)(\alpha,\beta):= 
\cX(\alpha,\beta)^{G\times G}$. 

\begin{remark} If $X\in\alpha$ and $Y\in \beta$, then $(\cC/G)(\alpha,\beta)
\simeq \cC(X,Y)^{G_{X}\times G_{Y}}$, where $G_{X} = \{g\in G|gX=X\}$. 
\end{remark}

\begin{proposition} \label{prop:factor}
If $\cD$ is a small category, then a functor $\cC/\Gamma\to\cD$
is the same as a functor $F:\cC\to\cD$, such that $F(gX)=F(X)$  and 
$F(i_{X}^{g}) = \on{id}_{F(X)}$ for any $g\in G$, 
$X\in\on{Ob}\cC$. 
\end{proposition}

The proof is immediate. 

\subsection{Quotients of the braid group} Let $B_{n}$ be the braid group of $n$
strands in the plane. It is presented by generators $\sigma_{1},\ldots,\sigma_{n-1}$
subject to the Artin relations $\sigma_{i}\sigma_{i+1}\sigma_{i} = \sigma_{i+1}
\sigma_{i}\sigma_{i+1}$ for $i=1,\ldots,n-2$ and $\sigma_{i}\sigma_{j} = 
\sigma_{j}\sigma_{i}$ for $|i-j|>1$. Its center $Z_n:=Z(B_{n})$ is isomorphic to $\ZZ$ and 
is generated by $(\sigma_{1}\cdots\sigma_{n-1})^{n}$. There is a morphism 
$B_{n}\to S_{n}$, uniquely determined by $\sigma_{i}\mapsto s_{i}:= (i,i+1)$; 
it factors through a morphism $B_{n}/Z_n\to S_{n}$. 

\begin{lemma}
Let $C_{n}:= \langle g|g^{n}=1\rangle$ be the cyclic group of order $n$. 
We have an injection $C_{n}\hookrightarrow S_{n}$ via $g\mapsto 
\bigl(\begin{smallmatrix} 1 & 2 & \ldots & n \\ 2 & 3 & \ldots 
& 1\end{smallmatrix}\bigl)$, which admits a lift $C_{n}\to B_{n}/Z_n$, 
given by $g\mapsto \sigma_{1}\cdots \sigma_{n-1}$. 
\end{lemma}

Let $\Gamma_{0,n}:= B_{n}/((\sigma_{1}\cdots\sigma_{n-1})^{n},\sigma_{1}\cdots
\sigma_{n-1}^{2}\cdots\sigma_{1})$ be the mapping class group of type $(0,n)$
(see \cite{Bi}). The relation $\sigma_{1}\cdots\sigma_{n-1}^{2}\cdots
\sigma_{1}=1$ is called the sphere relation as the quotient $B_{n}/(\sigma_{1}
\cdots \sigma_{n-1}^{2}\cdots\sigma_{1})$ is isomorphic to the braid group 
of $n$ points on the sphere. In this group, the relation 
$(\sigma_{1}\cdots\sigma_{n-1})^{2n}=1$ holds. The morphism 
$B_{n}\to S_{n}$ factors through a morphism $\Gamma_{0,n}\to S_{n}$. 

The dihedral group $D_{n}:= \langle r,s|r^{n}=s^{2}=(rs)^{2}=1\rangle$ may be viewed
as a subgroup of $S_{n}$ via $r\mapsto 
\bigl( \begin{smallmatrix} 1 & 2 & \ldots & n \\
2 & 3 & \ldots & 1\end{smallmatrix}\bigl)$, $s\mapsto 
\bigl( \begin{smallmatrix} 1 & 2 & \ldots & n \\
n & n-1 & \ldots & 1\end{smallmatrix}\bigl)$. 

\begin{lemma}
There exists a unique morphism $D_{n}\to \Gamma_{0,n}$, $r\mapsto \sigma_{1}
\cdots \sigma_{n-1}$, $s\mapsto \sigma_{1}(\sigma_{2}\sigma_{1})\cdots
(\sigma_{n-1}\cdots\sigma_{1})$, lifting the injection $D_{n}\hookrightarrow S_{n}$. 
\end{lemma}

{\em Proof.} One knows that $h_{n}:= \sigma_{1}(\sigma_{2}\sigma_{1})\cdots
(\sigma_{n-1}\cdots\sigma_{1})\in B_{n}$ is the half-twist, so that 
$h_{n}^{2} = z_{n} = (\sigma_{1}\cdots\sigma_{n-1})^{n} = \rho^{n}$, where 
$\rho = \sigma_{1}\cdots\sigma_{n-1}$ and $z_{n}$ is the full twist, 
generating $Z(B_{n})$. Moreover, $h_{n}\rho^{-1} = 
\on{im}(h_{n-1}\in B_{n-1}\to B_{n})$, so $(h_{n}\rho^{-1})^{2} = 
z_{n-1} = z_{n}(\sigma_{n-1}\cdots\sigma_{1}^{2}\cdots \sigma_{n-1})^{-1}$, 
where we identify $z_{n-1}$ with its image under $B_{n-1}\to B_{n}$. The images 
of $h_{n},\rho$ in $\Gamma_{0,n}$ therefore satisfy $\bar h_{n}^{2} = \bar\rho^{n}
=(\bar h_{n}\bar\rho^{-1})^{2}=1$, which are equivalent to a presentation of $D_{n}$. 
\hfill \qed\medskip 

\subsection{Teichm\"uller groupoids}

Let $G$ be a group and $\Gamma\subset S_{n}$ be a subgroup. Assume that $G\to S_{n}$
is a group morphism and let $\Gamma\to G$ be such that 
$\xymatrix{G\ar[r] & S_n \\ & \Gamma\ar[u]\ar[ul]}$
commutes. 
Let $S$ be a set, with $|S|=n$. 

Define a category $\cC_{G,S}$ by $\on{Ob}\cC_{G,S}:= \on{Bij}([n],S)$; 
for $\sigma,\sigma'\in \on{Ob}\cC_{G,S}$, $\on{Hom}(\sigma,\sigma'):= 
G\times_{S_{n}} \{(\sigma')^{-1}\sigma\}$; the composition of morphisms
is induced by the product in $G$. 

Define an action of $\Gamma$ on $\cC_{G,S}$ as follows. For $\gamma\in \Gamma$, 
$\sigma\in\on{Bij}([n],S)$, $\gamma \cdot \sigma:= \sigma\gamma^{-1}$, 
and $i_{\sigma}^{\gamma}\in \on{Hom}(\sigma,\sigma\gamma^{-1}) = 
G\times_{S_{n}}\{\gamma\}$ is $\on{im}(\gamma\in\Gamma\to G)$. 
We then obtain a quotient category $\cC_{\Gamma,G,S}:= \cC_{G,S}/\Gamma$. 

\begin{example} When $G = B_{n}/Z_n$ and $\Gamma = C_{n}$, we set 
$\ul{\on{Cyc}}(S):= \cC_{\Gamma,G,S}$; its set of objects is $\on{Cyc}(S):= 
\on{Bij}([n],S)/C_{n}$ (the set of cyclic orders on $S$). 
\end{example}

\begin{example} When $G = \Gamma_{0,n}$ and $\Gamma = D_{n}$, we set 
$\ul{\on{Dih}}(S):= \cC_{\Gamma,G,S}$; its set of objects is $\on{Dih}(S):= 
\on{Bij}([n],S)/D_{n} = \on{Cyc}(S)/\{\pm1\}$, which we call the set of dihedral orders
on $S$. 
\end{example}

\begin{definition} If $\cC$ is a small category and $T\stackrel{\pi}{\to}\on{Ob}\cC$
is a map, we define the category $\pi^{*}\cC$ by $\on{Ob}\pi^{*}\cC:= T$
and $\on\pi^{*}\cC(t,t'):= \cC(\pi(t),\pi(t'))$ for $t,t'\in T$. 
\end{definition}

We have natural maps 
$$
\{\text{planar 3-valent trees with leaves bijectively indexed by }S\}
\stackrel{\pi_{cyc}}{\to} \on{Cyc}(S)$$
and 
$$
\{\text{planar 3-valent trees with leaves bijectively indexed by }S\}/(\text{mirror 
symmetry}) \stackrel{\pi_{dih}}{\to} \on{Dih}(S). 
$$
We then set $T'_{0,S}:= \pi_{cyc}^{*}\ul{\on{Cyc}}(S)$, 
$T_{0,S}:=\pi_{dih}^{*}\ul{\on{Dih}}(S)$.  

When $S = [n]$, $T_{0,S}$ identifies with the fundamental groupoid to the moduli 
stack $M_{0,n}^{\QQ}$ with respect to the set of maximally degenerate real curves 
(see \cite{Sch}). 

\section{Contractions on (half-)balanced categories}

\subsection{(Half-)balanced categories}
Recall that a braided monoidal category (b.m.c.) is a set 
$(\cC,\otimes,{\bf 1},\beta_{XY},a_{XYZ})$, where $\cC$ is a category, 
$\otimes : \cC\times\cC\to\cC$ is a bifunctor, $\beta_{XY}:X\otimes Y\to 
Y\otimes X$ and $a_{XYZ}:(X\otimes Y)\otimes Z\to X\otimes (Y\otimes Z)$
are natural constraints, ${\bf 1}\in \on{Ob}\cC$ and $X\otimes {\bf 1}
\stackrel{\sim}{\to}X\stackrel{\sim}{\leftarrow}{\bf 1}\otimes X$ are natural 
isomorphisms, satisfying the hexagon and pentagon conditions (see e.g.~ \cite{Ka}). 

A balanced structure on the small 
b.m.c.~ $\cC$ is the datum of a natural assignment $\on{Ob}\cC
\ni X\mapsto \theta_{X}\in \on{Aut}_{\cC}(X)$, such that 
$$
\theta_{X\otimes Y} = (\theta_{X}\otimes\theta_{Y})\beta_{YX}\beta_{XY}
$$
for any $X,Y\in\on{Ob}\cC$ (see \cite{JS}); the naturality condition is $\theta_{X'}\phi 
= \phi\theta_{X}$ for any $X,X'\in\on{Ob}\cC$ and $\phi\in
\on{Hom}_{\cC}(X,X')$.

Similarly, a half-balanced structure on $\cC$ is the data of: (a) an involutive
autofunctor $*:\cC\to\cC$, $X\mapsto X^{*}$, such that 
$(X\otimes Y)^{*} = Y^{*}\otimes X^{*}$, $(f\otimes g)^{*} = g^{*}
\otimes f^{*}$ for any $X,\ldots,Y'\in\on{Ob}\cC$ and $f\in\on{Hom}_{\cC}(X,X')$, 
$g\in\on{Hom}_{\cC}(Y,Y')$, $a_{X}^{*} = a_{X^{*}}$, $\beta_{XY}^{*} = 
\beta_{Y^{*}X^{*}}$, $a_{XYZ}^{*}=a_{Z^{*}Y^{*}X^{*}}$; 
(b) a natural assignment $\on{Ob}\cC\in X\mapsto 
a_{X}\in \on{Iso}_{\cC}(X,X^{*})$, such that 
$$
a_{X\otimes Y} = (a_{Y}\otimes a_{X})\beta_{XY}$$
for any $X,Y\in\on{Ob}\cC$; here naturality means that 
$a_{Y}\phi = \phi^{*}a_{X}$ for any $\phi\in\on{Hom}_{\cC}(X,Y)$. 

Note that a half-balanced structure gives rise to a balanced structure by 
$\theta_{X}:= a_{X^{*}}a_{X}$. 

\subsection{Contractions}

\begin{definition} A contraction on the small balanced category $\cC$ is a functor 
$\langle -\rangle : \cC\to\cO$, $X\mapsto \langle X\rangle$, such that: 

1) for any $X,Y\in\on{Ob}\cC$, $\langle Y\otimes X\rangle 
= \langle X\otimes Y \rangle(=:\langle X,Y\rangle)$, and 
$\langle (\theta_{Y}\otimes\on{id}_{X})\beta_{XY}\rangle = 
\on{id}_{\langle X,Y\rangle}$; 

2) for any $X,Y,Z\in \on{Ob}\cC$, $\langle (X\otimes Y)\otimes Z
\rangle = \langle X\otimes (Y\otimes Z)\rangle (=:\langle X,Y,Z\rangle)$
 and $\langle a_{XYZ}\rangle= \on{id}_{\langle X,Y,Z\rangle}$. 
 \end{definition}

When needed, we will call such a contraction a ``balanced contraction''. 

\begin{remark} \label{rem:theta}
These axioms imply $\langle \theta_{X\otimes Y}\rangle 
= \on{id}_{\langle X,Y\rangle}$ for any $X,Y\in\on{Ob}\cC$, and therefore 
 $\langle \theta_{X}\rangle = \on{id}_{\langle X\rangle}$ by taking 
 $Y={\bf 1}$. 
\end{remark}

\begin{definition}
A contraction on the small half-balanced category $\cC$ is a functor $\langle-\rangle:
\cC\to\cO$, such that:

1) $\langle-\rangle$ is a balanced contraction on $\cC$; 

2) for any $X\in\on{Ob}\cC$, $\langle X\rangle = \langle X^{*}\rangle$
and $\langle a_{X}\rangle = \on{id}_{\langle X\rangle}$.  
\end{definition}

When needed, we will call such a contraction a ``half-balanced contraction''. 

\begin{lemma} \label{lemma:hbal}
If $\langle -\rangle : \cC\to \cO$ is a contraction on a half-balanced 
category, then for any $X,Y\in \on{Ob}\cC$, $\langle \theta_{X}\otimes 
\theta_{Y}^{-1}\rangle = \on{id}_{\langle X,Y
\rangle} = \langle (\theta_{X}^{2}\otimes \on{id}_{Y})
\beta_{YX}\beta_{XY}\rangle$. 
\end{lemma}

{\em Proof.}
We have 
\begin{eqnarray*}
& & \beta_{XY}^{-1}(\theta_{Y}^{-1}\otimes \on{id}_{X})
a_{X^{*}\otimes Y^{*}}\beta_{X^{*}Y^{*}}^{-1}(\theta_{Y^{*}}^{-1}
\otimes \on{id}_{X^{*}})a_{X\otimes Y}\\
& & = 
\beta_{XY}^{-1}(\theta_{Y}^{-1}\otimes \on{id}_{X})
a_{X^{*}\otimes Y^{*}}\beta_{X^{*}Y^{*}}^{-1} a_{X\otimes Y} 
(\on{id}_{X}\otimes \theta_{Y}^{-1})
\\
& & =
\beta_{XY}^{-1}(\theta_{Y}^{-1}\otimes \on{id}_{X})
(a_{Y^{*}}\otimes a_{X^{*}}) a_{X\otimes Y} 
(\on{id}_{X}\otimes \theta_{Y}^{-1})
\\
& & =
\beta_{XY}^{-1}(\theta_{Y}^{-1}\otimes \on{id}_{X})
(a_{Y^{*}}a_{Y}\otimes a_{X^{*}}a_{X}) \beta_{XY} 
(\on{id}_{X}\otimes \theta_{Y}^{-1})
\\
& & =
\beta_{XY}^{-1}(\on{id}_{Y}\otimes \theta_{X})\beta_{XY} 
(\on{id}_{X}\otimes \theta_{Y}^{-1}) = \theta_{X}\otimes \theta_{Y}^{-1}.  
\end{eqnarray*}
Now $\langle a_{X\otimes Y}\rangle = \langle a_{X^{*}\otimes Y^{*}}\rangle
= \on{id}_{\langle X,Y\rangle}$ by the half-balanced contraction axiom, and 
$\langle \beta_{XY}^{-1}(\theta_{Y}^{-1}\otimes \on{id}_{X})
\rangle = \langle \beta_{X^{*}Y^{*}}^{-1}(\theta_{Y^{*}}^{-1}
\otimes \on{id}_{X^{*}})\rangle = \on{id}_{\langle X,Y\rangle}$
by the balanced contraction axiom. It follows that $\langle
 \theta_{X}\otimes \theta_{Y}^{-1}\rangle = \on{id}_{\langle X,Y\rangle}$. 
The second statement follows from $ (\theta_{X}^{2}\otimes \on{id}_{Y})
\beta_{YX}\beta_{XY} = (\theta_{X}\otimes \theta_{Y}^{-1})\theta_{X
\otimes Y}$ and $\langle\theta_{X\otimes Y}\rangle = \on{id}_{\langle X,Y
\rangle}$. 
\hfill \qed\medskip 

\subsection{Relation with braid group representations}

Set $\tilde B_n:= \ZZ^n\rtimes B_n$, where the action of $B_n$
is $\ZZ^n$ is via $B_n \to S_n\to\on{Aut}(\ZZ^n)$; $\tilde B_n$
is usually called the framed braid group of the plane. If $\cC$ is 
a balanced b.m.c.~and $X\in\on{Ob}\cC$, then there is a morphism 
$\tilde B_n\to \on{Aut}_\cC(X^{\otimes n})$ (a parenthesization 
of the $n$th fold tensor product being chosen), given in the strict
case by 
$$
\delta_i\mapsto \on{id}_{X^{\otimes i-1}} \otimes\theta_X \otimes
\on{id}_{X^{\otimes n-i}}, \quad 
\sigma_i\mapsto \on{id}_{X^{\otimes i-1}} \otimes\beta_{X,X} \otimes
\on{id}_{X^{\otimes n-i-1}}. 
$$ 
Here $\delta_i$ is the $i$th generator of $\ZZ^n\subset \tilde B_n$. 

We now define $\widetilde{B_n/Z_n}$ to be the quotient of $\tilde B_n$
by its central subgroup (isomorphic to $\ZZ$) generated by 
$(\prod_{i=1}^n\delta_i)z_n$ (recall that $z_n$ is a generator of
$Z_n = Z(B_n)$; the product in $\ZZ^n$ is denoted multiplicatively). 
One can prove that there is a (generally non-split) exact sequence 
$1\to\ZZ^n\to \widetilde{B_n/Z_n}\to B_n/Z_n\to 1$. 

\begin{proposition}
Let $\cC\stackrel{\langle-\rangle}{\to}\cO$ be a balanced contraction
of $\cC$, then we have a commutative diagram 
$$
\begin{matrix}
B_n & \leftarrow & \tilde B_n & \to & \on{Aut}_\cC(X^{\otimes n}) \\
\downarrow &&\downarrow && \downarrow\scriptstyle{\langle-\rangle}\\
B_n/Z_n &\leftarrow & \widetilde{B_n/Z_n} &\to&
\on{Aut}_\cO(\langle X^{\otimes n}\rangle)\end{matrix}$$
\end{proposition}

{\em Proof.} We have $\on{im}((\prod_{i=1}^n\delta_i)z_n
\in \tilde B_n\to \on{Aut}_\cC(X^{\otimes n})) = \theta_{X^{\otimes n}}$, 
so according to Remark \ref{rem:theta}, the image of this in 
$\on{Aut}_\cO(\langle X^{\otimes n}\rangle)$ is 
$\on{id}_{\langle X^{\otimes n}\rangle}$. The factorization implied
in the right square follows. The left square obviously commutes. 
\hfill\qed\medskip 

Set now $\tilde\Gamma_{0,n}$ be the quotient of $\tilde B_n$ by the 
normal subgroup generated by $(\prod_{i=1}^n\delta_i)z_n$ and 
$\delta_1^2\sigma_1\cdots\sigma_{n-1}^2\cdots\sigma_1$. Then we have an 
exact sequence $1\to \ZZ^n\to\tilde\Gamma_{0,n}\to\Gamma_{0,n}\to 1$. 

\begin{proposition}
Let $\cC$ be a half-balanced b.m.c., let 
$\cC\stackrel{\langle-\rangle}{\to}\cO$ be a half-balanced contraction 
and let $X\in\on{Ob}\cC$. Then we have a comutative diagram 
$$
\begin{matrix}
B_n & \leftarrow & \tilde B_n & \to & \on{Aut}_\cC(X^{\otimes n}) \\
\downarrow &&\downarrow && \downarrow\scriptstyle{\langle-\rangle}\\
\Gamma_n &\leftarrow & \tilde\Gamma_n &\to&
\on{Aut}_\cO(\langle X^{\otimes n}\rangle)\end{matrix}
$$
\end{proposition}

{\em Proof.} We have
$\on{im}(\delta_1^2\sigma_1\cdots\sigma_{n-1}^2\cdots\sigma_1
\in \tilde B_n\to \on{Aut}_\cC(X^{\otimes n})) = (\theta_X^2\otimes
\on{id}_Y)\beta_{YX}\beta_{XY}$ ($Y=X^{\otimes n-1}$), whose image in 
$\on{Aut}_\cO(\langle X^{\otimes n}\rangle)$ is $\on{id}_{X^{\otimes n}}$
by Lemma \ref{lemma:hbal}. \hfill\qed\medskip 

\section{Universal (half-)balanced categories}

\subsection{Universal (strict) braided monoidal categories}

Recall that the small b.m.c.~$\cC$ is called strict iff $(X\otimes Y)\otimes Z = 
X\otimes (Y\otimes Z)(=X\otimes Y\otimes Z)$ and 
$a_{X,Y,Z} = \on{id}_{X\otimes Y\otimes Z}$ for any $X,Y,Z\in \on{Ob}\cC$. 
Following \cite{JS}, we associate a universal strict b.m.c.~${\bf B}_{S}$ 
to each set $S$. Its set of objects is 
$\on{Ob}{\bf B}_{S}:= \sqcup_{n\geq 0} S^{n}$; the tensor product is 
defined by $\ul s\otimes \ul s' = (s_{1},\ldots,s_{n},s'_{1},\ldots,s'_{n'})\in 
S^{n+n'}$ for $\ul s = (s_{1},\ldots,s_{n})\in S^{n}$, 
$\ul s' = (s'_{1},\ldots,s'_{n'})\in S^{n'}$. If $\ul s\in S^{n}$, $\ul s'\in S^{n'}$, 
then $\on{Hom}_{{\bf B}_{S}}(\ul s,\ul s') = \emptyset$ if $n\neq n'$, and 
$\on{Hom}_{{\bf B}_{s}}(\ul s,\ul s') = B_{n}\times_{S_{n}} \{f\in S_{n} | 
\ul s'f = \ul s\}$ if $n=n'$. The tensor product of morphisms is induced by restriction 
from the group morphism $B_{n}\times B_{n'}\to B_{n+n'}$, $(b,b')\mapsto b*b'$, 
uniquely determined by $\sigma_{i}*1 = \sigma_{i}$, $1*\sigma_{i'} = \sigma_{n-1+i'}$
(which corresponds to the juxtaposition of braids). The braiding is 
$\beta_{\ul s,\ul s'} = b_{nn'}$, where $b_{nn'}\in B_{n+n'}$ is 
given by 
$$
b_{nn'} = (\sigma_{n'}\cdots\sigma_{1})\cdots (\sigma_{n+n'-1}\cdots\sigma_{n}).
$$
The universal property of ${\bf B}_{S}$ is then expressed as follows: to each 
strict small b.m.c.~$\cC$ and any map $S\to \on{Ob}\cC$, there corresponds a unique 
tensor functor ${\bf B}_{S}\to \cC$, such that the diagram 
$\xymatrix{S \ar[r]\ar[dr]& \on{Ob}\cC \\ & \on{Ob}{\bf B}_S\ar[u]}$
commutes. 
 
We now describe the universal b.m.c.~${\bf PaB}_{S}$ associated to $S$
(\cite{JS,Ba}). Define first $T_{n}$ as the set of parenthesizations of a word in $n$
identical letters. Equivalently, this is the set of planar 3-valent rooted trees with $n$ leaves, 
e.g.~the tree 

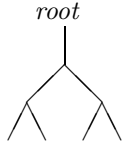
\setlength{\unitlength}{1mm}
\begin{picture}(60,20)(-30,0)
\put(21,16){{\it root}}
\put(25,10){\line(0,1){5}}
\put(25,10){\line(1,-1){5}}
\put(25,10){\line(-1,-1){5}}
\put(20,5){\line(1,-2){2.5}}
\put(20,5){\line(-1,-2){2.5}}
\put(30,5){\line(1,-2){2.5}}
\put(30,5){\line(-1,-2){2.5}}
\end{picture}

\noindent
corresponds to the word $(\bullet\bullet)(\bullet\bullet)$. 
The concatenation of words 
is a map $T_{n}\times T_{m}\to T_{n+m}$, $(t,t')\mapsto t*t'$ (e.g., 
$(\bullet\bullet,\bullet\bullet)\mapsto (\bullet\bullet)(\bullet)\bullet$); 
this is illustrated in terms of trees as follows 

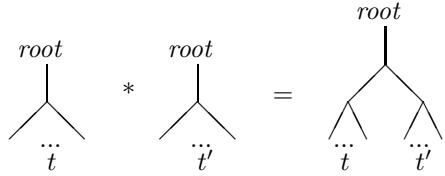
\setlength{\unitlength}{1mm}
\begin{picture}(60,25)(-30,0)
\put(1,16){{\it root}}
\put(5,10){\line(0,1){5}}
\put(5,10){\line(1,-1){5}}
\put(5,10){\line(-1,-1){5}}
\put(4,4){...}
\put(5,1){$t$}
\put(15,10){*}
\put(21,16){{\it root}}
\put(25,10){\line(0,1){5}}
\put(25,10){\line(1,-1){5}}
\put(25,10){\line(-1,-1){5}}
\put(24,4){...}
\put(25,1){$t'$}
\put(35,10){=}
\put(46,21){{\it root}}
\put(50,15){\line(0,1){5}}
\put(50,15){\line(1,-1){5}}
\put(50,15){\line(-1,-1){5}}
\put(45,10){\line(1,-2){2.5}}
\put(45,10){\line(-1,-2){2.5}}
\put(55,10){\line(1,-2){2.5}}
\put(55,10){\line(-1,-2){2.5}}
\put(43,4){...}
\put(53,4){...}
\put(44,1){$t$}
\put(54,1){$t'$}
\end{picture}

\noindent
The set of objects of ${\bf PaB}_{S}$ is then defined by $\on{Ob}{\bf PaB}_{S}:= 
\sqcup_{n\geq 0} T_{n}\times S^{n}$; the tensor product is defined by 
$(t,\ul s)\otimes (t',\ul s'):= (t*t',\ul s\otimes \ul s')$. The morphisms are defined by 
$\on{Hom}_{{\bf PaB}_{S}}((t,\ul s),(t',\ul s')):= \on{Hom}_{{\bf B}_{S}}(\ul s,
\ul s')$. The tensor product of morphisms and the braiding and associativity 
constraints are uniquely determined by the condition that the obvious functor 
${\bf PaB}_{S}\to{\bf B}_{S}$ is monoidal. In particular, $a_{XYZ}$ corresponds to 
$1\in B_{|X|+|Y|+|Z|}$, where $|(\ul s,t)| = n$ for $(\ul s,t)\in T_{n}\times S^{n}$. 
Then ${\bf PaB}_{S}$ has a universal property with respect to non-necessarily strict 
braided monoidal categories, analogous to that of ${\bf B}_{S}$. 

\subsection{Universal balanced categories} \label{sec:bal}

For $\ul s\in \on{Ob}{\bf B}_{S}$, set $\theta_{\ul s}:= z_{|\ul s|}
\in  \on{Aut}_{{\bf B}_{S}}(\ul s) \subset B_{|\ul s|}$. The assignment 
$\ul s\mapsto \theta_{\ul s}$ equips ${\bf B}_{S}$ with a balanced structure. 
We denote by ${\bf B}_{S}^{bal}$ the resulting balanced strict b.m.c. One
checks that it has the following universal property: 

\begin{lemma} \label{univ:B:bal}
To any balanced strict small b.m.c.~$\cC$ and any map $S\stackrel{f}{\to}
 \on{Ob}\cC$, such that $\theta_{f(s)} = \on{id}_{f(s)}$ for any $s\in S$, 
 there corresponds a unique functor ${\bf B}_{S}^{bal}\to \cC$ compatible with 
 the balanced and monoidal structures, such that the diagram 
 $\xymatrix{S \ar[r]\ar[dr]& \on{Ob}\cC \\ & \on{Ob}{\bf B}_{S}^{bal} \ar[u]}$
commutes. 
\end{lemma}

If now $X = (\ul s,t)\in \on{Ob}{\bf PaB}_{S}$, we set $\theta_{X}:= 
\theta_{\ul s}\in \on{Aut}_{{\bf B_{S}}}(\ul s) = \on{Aut}_{{\bf PaB}_{S}}(X)$. 
The assignment $X\mapsto \theta_{X}$ equips ${\bf PaB}_{S}$ with the 
structure of a balanced b.m.c., denoted ${\bf PaB}_{S}^{bal}$ and enjoying a universal 
property with respect to maps $S\to \on{Ob}\cC$, where 
$\cC$ is a balanced braided monoidal category such that 
$\theta_{f(s)} = \on{id}_{f(s)}$ similar to Lemma \ref{univ:B:bal}. 

\subsection{Universal half-balanced categories}

We define an involution $*:{\bf B}_{S}\to {\bf B}_{S}$ as follows. 
It is given at the level of objects by $\ul s^{*}:= (s_{n},\ldots,s_{1})$
for $\ul s = (s_{1},\ldots,s_{n})$ and the level of morphisms by 
restriction of the automorphism $\sigma_{i}\mapsto \sigma_{n-i}$
of $B_{n}$. For $\ul s\in \on{Ob}{\bf B}_{S}$, we set $a_{\ul s}:= 
h_{|\ul s|}\in  \on{Iso}_{{\bf B}_{S}}(\ul s,\ul s^{*})
\subset B_{|\ul s|}$. This defines the structure of a half-balanced 
category on ${\bf B}_{S}$, denoted ${\bf B}_{S}^{hbal}$, whose 
balanced structure is that described in Subsection \ref{sec:bal}. It has 
the following universal property: 

\begin{lemma} \label{univ:B:hbal}
For each strict half-balanced small b.m.c.~$\cC$ and each 
map $S\stackrel{f}{\to} \on{Ob}\cC$ such that for any $s\in S$, 
$f(s)^{*} = f(s)$ and $a_{f(s)} = \on{id}_{f(s)}$, there exists a unique 
functor ${\bf B}_{S}^{hbal}\to \cC$, compatible with the 
monoidal and half-balanced structures, and such that the diagram 
 $\xymatrix{S \ar[r]\ar[dr]& \on{Ob}\cC 
\\ & \on{Ob}{\bf B}_{S}^{bal} \ar[u]}$
commutes. 
\end{lemma}

We now define an involution $*$ of ${\bf PaB}_{S}$ as follows. At the level of 
objects, it is given by $X^{*}=(t^{*},\ul s^{*})$ for $X=(t,\ul s)$, where 
$t^{*}$ is the parenthesized word $t$, read in the reverse order (in terms of trees, 
this is the mirror image of $t$). At the level of
morphisms, it coincides with the involution $*$ of ${\bf B}_{S}$. We define 
the assignment $\on{Ob}{\bf PaB}_{S}\ni X\mapsto a_{X}$ by 
$a_{X}:=a_{\ul s}\in \on{Iso}_{{\bf B}_{S}}(\ul s,\ul s ^{*})=
\on{Iso}_{{\bf PaB}_{S}}(X,X^{*})$ for $X = (t,\ul s)$. 
This equips ${\bf PaB}_{S}$ with a half-balanced structure; the resulting 
half-balanced b.m.c.~is denoted ${\bf PaB}_{S}^{hbal}$. Its underlying 
balanced b.m.c.~is ${\bf PaB}_{S}^{bal}$. It has a universal property 
with respect to half-balanced small braided monoidal categories $\cC$ and maps 
$S\stackrel{f}{\to}\on{Ob}\cC$, such that $f(s)^{*}=f(s)$ and 
$a_{f(s)} = \on{id}_{f(s)}$, similar to that of Lemmas \ref{univ:B:bal} and 
\ref{univ:B:hbal}. 

\section{Universal contractions for balanced categories}

We will construct categories ${\bf (Pa)Cyc}_{S}$ and a diagram 
$\begin{matrix}{\bf PaB}_{S}^{bal} & \to & {\bf PaCyc}_{S} \\
\downarrow && \downarrow \\
{\bf B}_{S}^{bal} & \to & {\bf Cyc}_{S}\end{matrix}$ in which the 
horizontal functors are contractions and the left vertical functor is the canonical 
monoidal functor. 

We construct ${\bf Cyc}_{S}$ as follows. Define first $\widetilde{\bf Cyc}_{S}$
as the category with the same objects as ${\bf B}_{S}^{bal}$, and $B_{n}$ replaced by 
$B_{n}/Z_n)$ in the definition of morphisms. Define an action of $\ZZ$ on 
$\widetilde{\bf Cyc}_{S}$ by $1\cdot (s_{1},\cdots,s_{n}):= (s_{n},s_{1},\ldots,s_{n-1})$ and 
$i_{\ul s}^{1}\in \on{Iso}(\ul s,1\cdot \ul s) \subset B_{n}/Z_n$ 
is the class of $\sigma_{1}\cdots\sigma_{n-1}$. 
We then set ${\bf Cyc}_{S}:= \widetilde{\bf Cyc}_{S}/\ZZ$. Note that 
$\on{Ob}{\bf Cyc}_{S} = \sqcup_{n\geq 0}\on{Cyc}_{n}(S)$, where
$\on{Cyc}_{n}(S) = S^{n}/C_{n}$. We then define a functor 
${\bf B}_{S}^{bal}\to{\bf Cyc}_{S}$ as the composite functor 
${\bf B}_{S}^{bal}\to \widetilde{\bf Cyc}_{S}\to{\bf Cyc}_{S}$. 

Let us show that the functor $\langle-\rangle:
{\bf B}_{S}^{bal}\to{\bf Cyc}_{S}$ satisfies the 
balanced contraction condition. If $\ul s,\ul s'\in\on{Ob}{\bf B}_{S}$, 
with $\ul s = (s_{1},\ldots,s_{n})$ and $\ul s' = (s'_{1},\ldots,s'_{n'})$, 
then $\ul s'\otimes \ul s = (s'_{1},\ldots,s_{n}) = n'\cdot (\ul s\otimes \ul s')$, 
which implies that $\langle \ul s \otimes \ul s'\rangle =
\langle \ul s' \otimes \ul s\rangle$. Then $(\theta_{\ul s'}\otimes
\on{id}_{\ul s})\beta_{\ul s,\ul s'}\in \on{Iso}_{{\bf B}_{S}^{bal}}
(\ul s \otimes \ul s',\ul s'\otimes \ul s) = B_{n+n'}$ corresponds to 
$(z_{n'}*\on{id}_{n})b_{nn'} = (\sigma_{1}\cdots \sigma_{n+n'-1})^{n'}$. 
Its image in $\widetilde{\bf Cyc}_{S}$ is then $i^{n'}_{\ul s \otimes \ul s'}
\in \widetilde{\bf Cyc}_{S}(\ul s \otimes \ul s',n'\cdot(\ul s \otimes \ul s'))$, 
whose image in ${\bf Cyc}_{S}$ is $\on{id}_{\langle \ul s,\ul s'\rangle}$. 

We now prove the universality of this contraction. 

\begin{proposition} \label{prop:univ:tr}
Let $\cC$ be a strict small balanced b.m.c., equipped with a map 
$S\stackrel{f}{\to}\on{Ob}\cC$ and a balanced contraction $\cC\to\cO$. 
Then there is a functor ${\bf Cyc}_{S}\to \cO$, such that the diagram 
$\begin{matrix} {\bf B}_{S}^{bal} & \to & {\bf Cyc}_{S} \\
\downarrow &&\downarrow \\
\cC & \to & \cO\end{matrix}$ commutes. 
\end{proposition}

{\em Proof.} First note that since $\langle \theta_{X}\rangle 
= \on{id}_{\langle X\rangle}$ for $X = f(s_{1})\otimes \cdots 
\otimes f(s_{n})$ and any $(s_{1},\ldots,s_{n})\in \on{Ob}{\bf B}_{S}^{bal}$, 
we have a functor $F : \widetilde{{\bf Cyc}}_{S}\to \cO$, such that the 
diagram $\begin{matrix} {\bf B}_{S}^{bal} & \to & \widetilde{\bf Cyc}_{S} \\
\downarrow &&\downarrow \\
\cC & \to & \cO\end{matrix}$ commutes. 

If $(s_{1},\ldots,s_{n})\in \on{Ob}\widetilde{\bf Cyc}_{S} 
= \on{Ob}{\bf B}_{S}^{bal}$, then $F(s_{1},\ldots,s_{n}) = \langle
f(s_{1})\otimes \cdots \otimes f(s_{n}) \rangle = \langle f(s_{n})
\otimes f(s_{1})\otimes \cdots \otimes f(s_{n-1})\rangle = F(s_{n},\ldots,
s_{n-1})$, therefore $F(gX)=F(X)$ for any $X\in\on{Ob}\widetilde{\bf Cyc}_{S}$
and any $g\in\ZZ$. Moreover, we have 
\begin{eqnarray*}
&& F(i^{1}_{(s_{1},\ldots,s_{n})}) = F(\sigma_{1}\cdots \sigma_{n-1}) = 
\langle (\theta_{f(s_{n})}\otimes \on{id}_{f(s_{1})\otimes 
\cdots \otimes f(s_{n-1})})\beta_{f(s_{1})\otimes \cdots \otimes f(s_{n-1}),
f(s_{n})}\rangle \\ && = \on{id}_{F(s_{1},\cdots,s_{n})}
\end{eqnarray*} by the balanced contraction 
property. 

According to Proposition \ref{prop:factor}, this implies that we have a 
factorization $\xymatrix{\widetilde{\bf Cyc}_{S} \ar[r]\ar[dr] & {\bf Cyc}_{S}
\ar[d] \\ & \cO}$
\hfill \qed\medskip

We now construct the category ${\bf PaCyc}_{S}$ as follows. Let 
$PlT_{n}:= \{$planar 3-valent trees equipped with a bijection $\{$leaves$\}\to
[n]$, compatible with the cyclic orders$\}$. We first define the category 
$\widetilde{\bf PaCyc}_{S}$ by $\on{Ob}\widetilde{\bf PaCyc}_{S} = 
\sqcup_{n\geq 0} PlT_{n}\times S^{n}$, $\on{Hom}_{\widetilde{\bf PaCyc}_{S}}
((t,\sigma),(t',\sigma')) = \on{Hom}_{\widetilde{{\bf Cyc}}_{S}}(\sigma,\sigma')$. 
We define an action of $\ZZ$ on ${\bf PaCyc}_{S}$ by $1\cdot(t,(s_{1},\ldots,s_{n}))
:= (t',(s_{n},s_{1},\ldots,s_{n-1}))$, where if $t = (T,\{$leaves of 
$T\}\stackrel{\alpha}{\to}[n])$, then $t':= (T,\{$leaves of $T\}
\stackrel{\alpha}{\to}[n]\stackrel{+1\on{\ mod\ }n}{\to}[n])$, 
and $i^{1}_{(t,\sigma)}:= i^{1}_{\sigma}$; we then set ${\bf PaCyc}_{S}:= 
\widetilde{{\bf PaCyc}}_{S}/\ZZ$, so in particular $\on{Ob}{\bf PaCyc}_{S}
= \{($a planar 3-valent tree, a map $\{$leaves$\}\to S)\}$.

We define a map $T_{n}\to PlT_{n}$, $t\mapsto \pi(t)$ as the operation of (a) 
assigning labels $1,\ldots,n$ to the vertices of the tree $t$, numbered from left to right; 
(b) replacing the root and the edges connected to it, by a single edge. E.g., we
have 

\setlength{\unitlength}{1mm}
\begin{picture}(60,20)(-30,0)
\put(0,7){$\pi\Big($}\put(11,16){{\it root}}
\put(15,10){\line(0,1){5}}
\put(15,10){\line(1,-1){5}}
\put(15,10){\line(-1,-1){5}}
\put(10,5){\line(1,-2){2.5}}
\put(10,5){\line(-1,-2){2.5}}
\put(20,5){\line(1,-2){2.5}}
\put(20,5){\line(-1,-2){2.5}}
\put(25,7){$\Big) = $}
\put(35,10){\line(1,0){10}}
\put(35,10){\line(1,-2){2.5}}
\put(35,10){\line(-1,-2){2.5}}
\put(45,10){\line(1,-2){2.5}}
\put(45,10){\line(-1,-2){2.5}}
\put(32,2){1}\put(42,2){3}
\put(37,2){2}\put(47,2){4}
\end{picture}

We define a functor ${\bf PaB}_{S}^{bal}\to {\bf PaCyc}_{S}$ by the condition 
that (a) at the level of objects, it is given by the map $\sqcup_{n\geq 0}
T_{n}\times S^{n}\to \sqcup_{n\geq 0} (PlT_{n}\times S^{n})/C_{n}$
and by projection, and (b) the diagram $\begin{matrix}
{\bf PaB}_{S}^{bal} &\to  &{\bf PaCyc}_{S} \\
\downarrow &  &\downarrow \\
{\bf B}_{S}^{bal}&\to &{\bf Cyc}_{S} \end{matrix}$ commutes. 
Let us check that this defines a contraction. $\langle X\otimes Y\rangle
= \langle Y\otimes X\rangle$ follows from the fact that for 
$t\in T_{n}$, $t'\in T_{n'}$, $\pi(t\otimes t')$ and $\pi(t'\otimes t)$
can be obtained from one another by cyclic permutation of $[n+n']$; 
here we recall that $(t,t')\mapsto t*t'$ is the concatenation map
$T_{n}\times T_{n'}\to T_{n+n'}$. The fact that $\langle (X\otimes Y)\otimes 
Z\rangle = \langle X\otimes (Y\otimes Z)\rangle$ follows from
$\pi((t*t')*t'') = \pi(t*(t'*t''))$, which  is illustrated as follows

\setlength{\unitlength}{1mm}
\begin{picture}(60,25)(-30,0)
\put(0,15){\line(1,-1){5}}
\put(0,15){\line(-1,-1){5}}
\put(0,15){\line(3,-1){15}}
\put(-5,10){\line(1,-2){2.5}}
\put(-5,10){\line(-1,-2){2.5}}
\put(-7,4){...}
\put(-6,1){$t$}
\put(5,10){\line(1,-2){2.5}}
\put(5,10){\line(-1,-2){2.5}}
\put(4,4){...}
\put(4,1){$t'$}
\put(15,10){\line(1,-2){2.5}}
\put(15,10){\line(-1,-2){2.5}}
\put(14,4){...}
\put(14,1){$t''$}
\put(22,9){$=$}
\put(45,15){\line(1,-1){5}}
\put(45,15){\line(-1,-1){5}}
\put(45,15){\line(-3,-1){15}}
\put(30,10){\line(1,-2){2.5}}
\put(30,10){\line(-1,-2){2.5}}
\put(28,4){...}
\put(29,1){$t$}
\put(40,10){\line(1,-2){2.5}}
\put(40,10){\line(-1,-2){2.5}}
\put(39,4){...}
\put(39,1){$t'$}
\put(50,10){\line(1,-2){2.5}}
\put(50,10){\line(-1,-2){2.5}}
\put(49,4){...}
\put(49,1){$t''$}
\put(57,9){$=$}
\put(75,15){\line(0,-1){5}}
\put(75,15){\line(2,1){5}}
\put(75,15){\line(-2,1){5}}
\put(75,10){\line(1,-2){2.5}}
\put(75,10){\line(-1,-2){2.5}}
\put(73,4){...}
\put(74,1){$t$}
\put(80,17.5){\line(1,0){5}}
\put(80,17.5){\line(1,2){2.5}}
\put(70,17.5){\line(-1,0){5}}
\put(70,17.5){\line(-1,2){2.5}}
\put(84,21.5){$\cdot$}
\put(85,20){$\cdot$}
\put(86,18.5){$\cdot$}
\put(87,22){$t'$}
\put(63.5,18.5){$\cdot$}\put(64.5,20){$\cdot$}\put(65.5,21.5){$\cdot$}
\put(62,22){$t''$}
\end{picture}

It is then clear that $\langle a_{XYZ}\rangle = \on{id}_{\langle X,Y,Z\rangle}$. 
The proof of $\langle (\theta_{Y}\otimes \on{id}_{X})\beta_{XY}\rangle = 
\on{id}_{\langle X,Y\rangle}$ is as above. 
We now prove the universality of the contraction $\langle-\rangle : 
{\bf PaB}_{S}^{bal}\to {\bf PaCyc}_{S}$. 

\begin{proposition} \label{bmc}
Let $\cC$ be a balanced small b.m.c., equipped with a contraction $\cC\to\cO$ and a 
map $S\to \on{Ob}\cC$. Then there exists a functor ${\bf PaCyc}_{S}\to\cO$, 
such that the diagram $\begin{matrix} {\bf PaB}_{S}^{bal} & \to & 
{\bf PaCyc}_{S} \\ \downarrow & & \downarrow \\
\cC & \to & \cO
\end{matrix}$ commutes. 
\end{proposition}

{\em Proof.} We first construct a functor $\widetilde{\bf PaCyc}_{S}\to\cO$, 
such that $\begin{matrix} {\bf PaB}_{S}^{bal} & \to & 
\widetilde{\bf PaCyc}_{S} \\ \downarrow & & \downarrow \\
\cC & \to & \cO \end{matrix}$ commutes. We define a map 
$PlT_{n}\times S^{n}\to \on{Ob}\cO$ as follows. Let $(t,(s_{1},\ldots,s_{n}))
\in PlT_{n}\times S^{n}$. Let $e$ be an edge of $t$. Cutting $t$ at $e$, we obtain 
two rooted trees $t_{i}$ ($i=1,2$) equipped with injective maps $\{$leaves of $t_{i}\}
\to [n]$. The images of these maps are of the form $\{a,a+1,\ldots,a+n_{1}\}$
and $\{a+n_{1}+1,\ldots,a+n_{1}+n_{2}\}$ (the integers being taken modulo 
$n$). We then define the image of $(t,(s_{1},\ldots,s_{n}))$ to be 
$\langle (\otimes_{i\in a+[n_{1}]}^{t_{1}}f(s_{i}))\otimes
(\otimes_{i\in a+n_{1}+[n_{2}]}^{t_{2}}f(s_{i}))\rangle$. The axioms 
then imply that this object do not depend on $e$. Indeed, if $e'$ is another 
edge, then to the shortest path $e = e_{1}\to e_{2}\to \cdots \to e_{k} = e'$
from $e$ to $e'$ there corresponds a sequence of isomorphisms of the 
corresponding objects; each isomorphism has the form 
$\langle A \otimes (B\otimes C)\rangle \stackrel{\langle a^{-1}_{ABC}
\rangle}{\longrightarrow} \langle (A\otimes B) \otimes C \rangle 
\stackrel{\langle\beta_{C,A\otimes B}^{-1}(\theta_{C}^{-1}\otimes
\on{id}_{A\otimes B})\rangle}{-\!\!\!-\!\!\!-\!\!\!-\!\!\!-\!\!\!-\!\!\longrightarrow} 
\langle C\otimes (A\otimes B)\rangle$, see

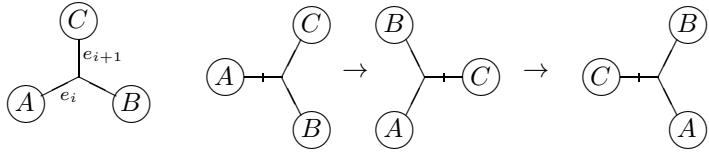
\setlength{\unitlength}{1mm}
\begin{picture}(60,30)(-30,0)
\put(-15,22.5){\circle{5}}
\put(-16.5,21.5){$C$}
\put(-15,15){\line(0,1){5}}
\put(-15,15){\line(2,-1){5}}
\put(-15,15){\line(-2,-1){5}}
\put(-22,11.5){\circle{5}}
\put(-23.5,10.5){$A$}
\put(-8,11.5){\circle{5}}
\put(-9.5,10.5){$B$}
\put(-14.5,17.5){$\scriptstyle{e_{i+1}}$}
\put(-17.5,12){$\scriptstyle{e_i}$}
\put(12,15){\line(-1,0){5}}
\put(4.5,15){\circle{5}}
\put(3,14){$A$}
\put(9.5,14.5){\line(0,1){1}}
\put(12,15){\line(1,2){2.5}}
\put(12,15){\line(1,-2){2.5}}
\put(16,22){\circle{5}}
\put(16,8){\circle{5}}
\put(14.5,21){$C$}
\put(14.5,7){$B$}
\put(20,15){$\to$}
\put(25.5,7){$A$}
\put(25.5,21){$B$}
\put(27,8){\circle{5}}
\put(27,22){\circle{5}}
\put(31,15){\line(-1,-2){2.5}}
\put(31,15){\line(-1,2){2.5}}
\put(31,15){\line(1,0){5}}
\put(33.5,14,5){\line(0,1){1}}
\put(38.5,15){\circle{5}}
\put(37,14){$C$}
\put(44,15){$\to$}
\put(62,15){\line(-1,0){5}}
\put(54.5,15){\circle{5}}
\put(53,14){$C$}
\put(59.5,14.5){\line(0,1){1}}
\put(62,15){\line(1,2){2.5}}
\put(62,15){\line(1,-2){2.5}}
\put(66,22){\circle{5}}
\put(66,8){\circle{5}}
\put(64.5,21){$B$}
\put(64.5,7){$A$}
\end{picture}

or $\langle A \otimes (B\otimes C)\rangle
\to \langle C\otimes (A\otimes B)\rangle\to 
\langle B\otimes (C\otimes A)\rangle$, see

\setlength{\unitlength}{1mm}
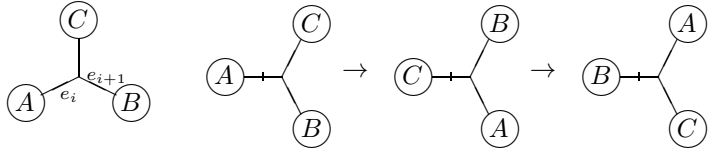
\begin{picture}(60,30)(-30,0)
\put(-15,22.5){\circle{5}}
\put(-16.5,21.5){$C$}
\put(-15,15){\line(0,1){5}}
\put(-15,15){\line(2,-1){5}}
\put(-15,15){\line(-2,-1){5}}
\put(-22,11.5){\circle{5}}
\put(-23.5,10.5){$A$}
\put(-8,11.5){\circle{5}}
\put(-9.5,10.5){$B$}
\put(-14,14.5){$\scriptstyle{e_{i+1}}$}
\put(-17.5,12){$\scriptstyle{e_i}$}
\put(12,15){\line(-1,0){5}}
\put(4.5,15){\circle{5}}
\put(3,14){$A$}
\put(9.5,14.5){\line(0,1){1}}
\put(12,15){\line(1,2){2.5}}
\put(12,15){\line(1,-2){2.5}}
\put(16,22){\circle{5}}
\put(16,8){\circle{5}}
\put(14.5,21){$C$}
\put(14.5,7){$B$}
\put(20,15){$\to$}
\put(37,15){\line(-1,0){5}}
\put(29.5,15){\circle{5}}
\put(28,14){$C$}
\put(34.5,14.5){\line(0,1){1}}
\put(37,15){\line(1,2){2.5}}
\put(37,15){\line(1,-2){2.5}}
\put(41,22){\circle{5}}
\put(41,8){\circle{5}}
\put(39.5,21){$B$}
\put(39.5,7){$A$}
\put(45,15){$\to$}
\put(62,15){\line(-1,0){5}}
\put(54.5,15){\circle{5}}
\put(53,14){$B$}
\put(59.5,14.5){\line(0,1){1}}
\put(62,15){\line(1,2){2.5}}
\put(62,15){\line(1,-2){2.5}}
\put(66,22){\circle{5}}
\put(66,8){\circle{5}}
\put(64.5,21){$A$}
\put(64.5,7){$C$}
\end{picture}
 
One then proves as before that we have a functor $\widetilde{\bf PaCyc}_{S}
\to \cO$, which factors through the action of $\ZZ$. \hfill \qed\medskip 
 
\section{Universal contractions for half-balanced categories}

We now construct categories ${\bf (Pa)Dih}_{S}$ and a commutative diagram  
$\begin{matrix} {\bf PaB}_{S}^{hbal} & \to & {\bf PaDih}_{S} \\
\downarrow  & & \downarrow \\
{\bf B}_{S}^{hbal} & \to & {\bf Dih}_{S}\end{matrix}$
where the horizontal functors are contractions. 

We first construct ${\bf Dih}_{S}$ as follows. Define first $\widetilde{\bf Dih}_{S}$
as the category with the same objects as ${\bf B}_{S}^{hbal}$, with $B_{n}$
replaced by its quotient $\Gamma_{0,n}$. Let $D:= \ZZ\rtimes (\ZZ/2)$ be the 
infinite dihedral group presented as $D:= \langle r,s|s^{2}=(rs)^{2}=1\rangle$. 
We define an action of $D$ on $\widetilde{\bf Dih}_{S}$ as follows. The action on 
objects is defined by $r\cdot (s_{1},\ldots,s_{n}):= (s_{n},s_{1},\ldots,s_{n-1})$, 
$s\cdot(s_{1},\ldots,s_{n}):=(s_{n},\ldots,s_{1})$, and $i_{\ul s}^{r} = \sigma_{1}
\cdots \sigma_{n-1}$, $i_{\ul s}^{s} = h_{n}$. We then set 
${\bf Dih}_{S} = \widetilde{\bf Dih}_{S}/D$.   
 
Note that $\on{Ob}{\bf Dih}_{S} = \sqcup_{n\geq 0}\on{Dih}_{n}(S)$, 
where $\on{Dih}_{n}(S) = S^{n}/D_{n}$, and $D_{n}$ is the 
quotient of $D$ by the relation $r^{n}=1$. We define a functor ${\bf B}_{S}^{hbal}
\stackrel{\langle-\rangle}{\to}{\bf Dih}_{S}$ as the composite functor 
${\bf B}_{S}^{hbal}\to \widetilde{\bf Dih}_{S}\to {\bf Dih}_{S}$. 
Let us show that it satisfies the half-balanced contraction conditions. 

We have a commutative diagram $\begin{matrix}{\bf B}_{S}^{bal} & \stackrel{
\langle -\rangle}{\to} & {\bf Cyc}_{S} \\
\downarrow & & \downarrow \\
{\bf B}_{S}^{hbal} & \stackrel{\langle-\rangle}{\to} & 
{\bf Dih}_{S}\end{matrix}$  
Since the left vertical functor is surjective on objects and the bottom functor  
is a balanced contraction, the upper functor satisfies the balanced contraction condition. 
If now $\ul s = (s_{1},\ldots,s_{n})\in \on{Ob}{\bf B}_{S}^{hbal}$, then 
$\ul s^{*} = (s_{n},\ldots,s_{1}) = s\cdot \ul s$, so the classes of $\ul s$
and $\ul s^{*}$ are the same in ${\bf Dih}_{S} = \widetilde{\bf Dih}_{S}/D$, 
hence $\langle \ul s\rangle = \langle \ul s^{*}\rangle$. 
Then $\langle a_{\ul s}\rangle = \langle i_{\ul s}^{s}\rangle = \on{id}_{
\langle \ul s\rangle}$. All this shows that ${\bf B}_{S}^{hbal}\stackrel{
\langle-\rangle}{\to}{\bf Dih}_{S}$ is a half-balanced contraction. We now prove 
the universality of this contraction. 

\begin{proposition} \label{prop:hrbmc}
Let $\cC$ be a strict half-balanced b.m.c., equipped with a map $S\stackrel{f}{\to}
\on{Ob}\cC$, such that $f(s)^{*} = f(s)$ for any $s\in S$, and with a half-balanced 
contraction $\cC\to\cO$. Then there exists a functor ${\bf Dih}_{S}\to \cO$, such that 
the diagram $\begin{matrix}{\bf B}_{S}^{hbal} & \to & {\bf Dih}_{S} \\
\downarrow & & \downarrow \\
\cC & \to & \cO\end{matrix}$ commutes. 
\end{proposition}

{\em Proof.} We define a functor $\widetilde{\bf Dih}_{S}\to\cO$ by the
following conditions: it coincides at the level of objects with the functor 
${\bf B}_{S}^{hbal}\to \cC\to\cO$; since the images by this functor of 
$z_{n},\sigma_{1}\cdots\sigma_{n-1}^{2}\cdots\sigma
\in \on{Aut}_{{\bf B}_{S}^{hbal}}(s_{1},\ldots,s_{n})
\subset B_{n}$ are respectively $\langle \theta_{f(s_{1})\otimes
\cdots \otimes f(s_{n})\rangle}$ and 
$$
\langle (\theta^{2}_{f(s_{1})}\otimes \on{id}_{\otimes_{i=2}^{n}f(s_{i})}
\beta_{\otimes_{i=2}^{n}f(s_{i}),f(s_{1})}
\beta_{f(s_{1}),\otimes_{i=2}^{n}f(s_{i})}
\rangle\in \on{Aut}_{\cO}(\langle f(s_{1})\otimes \cdots 
\otimes f(s_{n})\rangle), 
$$
which are the identity by Remark \ref{rem:theta} and Lemma \ref{lemma:hbal}, 
the composite functor ${\bf B}_{S}^{hbal}\to \cC\to\cO$ factorizes as 
$\begin{matrix}{\bf B}_{S}^{hbal} &\to&\widetilde{\bf Dih}_{S} \\
\scriptstyle{f}\downarrow & & \downarrow \\
\cC &\to &\cO\end{matrix}$ 
We now show as above that $F$ factorizes as 
$\begin{matrix}\widetilde{\bf Dih}_{S} &\to& {\bf Dih}_{S} \\
 & \scriptstyle{F}\searrow& \downarrow \\
 & &\cO\end{matrix}$ 
Indeed, for $\ul s = (s_{1},\ldots,s_{n})\in \on{Ob}\widetilde{\bf Dih}_{S}$, 
then $F(\ul s) = \langle f(s_{1})\otimes \cdots \otimes f(s_{n})\rangle$. 
Then $F(r\cdot \ul s) = \langle f(s_{n})\otimes \cdots \otimes f(s_{n-1})
\rangle = F(\ul s)$, using the axiom $\langle X \otimes Y\rangle = 
\langle Y\otimes X\rangle$ of balanced contractions, and $F(s\cdot \ul s) = 
\langle f(s_{n})\otimes \cdots \otimes f(s_{1})\rangle = \langle 
(f(s_{1})\otimes \cdots \otimes f(s_{n}))^{*}\rangle = 
\langle f(s_{1})\otimes \cdots \otimes f(s_{n})\rangle = F(\ul s)$
using the axiom $\langle X^{*}\rangle = \langle X\rangle$ of half-balanced 
contraction. If now $\ul s = (s_{1},\cdots,s_{n})\in \on{Ob}\widetilde{\bf Dih}_{S}$, 
then $F(i_{\ul s}^{r}) = \on{id}_{\langle X\rangle}$ by the same argument 
as in Proposition \ref{prop:univ:tr}, and 
\begin{eqnarray*} F(i_{\ul s}^{s}) = f(h_{n}) = 
&& (a_{f(s_{n})}\otimes \cdots \otimes a_{f(s_{1})})
(\on{id}_{f(s_{1})\otimes \cdots \otimes f(s_{n-2})}\otimes 
\beta_{f(s_{n-1}),f(s_{n})})
\\
 && (\on{id}_{f(s_{1})\otimes \cdots \otimes f(s_{n-3})}\otimes 
\beta_{f(s_{n-2}),f(s_{n-1})\otimes f(s_{n})})\cdots 
\beta_{f(s_{1}),f(s_{2})\otimes\cdots\otimes f(s_{n})} \\
 && = a_{f(s_{1})
\otimes \cdots \otimes f(s_{n})}.
\end{eqnarray*} 
Hence $F(i_{\ul s}^{s}) = \langle a_{f(s_{1})
\otimes \cdots \otimes f(s_{n})}\rangle = \on{id}_{\langle \ul s \rangle}$. 
So we have the desired factorization of $F$. \hfill \qed\medskip 

We now construct the category ${\bf PaDih}_{S}$ as follows. We first define the 
category $\widetilde{\bf PaDih}_{S}$ by $\on{Ob}\widetilde{\bf PaDih}_{S}
= \sqcup_{n\geq 0} PlT_{n}\times S^{n}$, $\widetilde{\bf PaDih}_{S}
((t,\ul s),(t',\ul s')) = \widetilde{\bf Dih}_{S}(\ul s,\ul s')$. The group $D$ acts on 
$\widetilde{\bf PaDih}_{S}$ as follows. The action on objects is $g\cdot (t,\ul s)
= (g\cdot t,g\cdot \ul s)$, where for $t = (T,\{$leaves of $T\}
\stackrel{\alpha}{\to}[n])$, $r\cdot t = (T,\{$leaves of $T\}
\stackrel{\alpha}{\to}[n]\stackrel{+1 \on{\ mod\ } n}{\to}[n])$, 
 $s\cdot t = (T,\{$leaves of $T\}
\stackrel{\alpha}{\to}[n]\stackrel{x\mapsto n+1-x}{\to}[n])$, 
and $i_{(t,\ul s)}^{g} = i_{\ul s}^{g}\in \on{Iso}_{\widetilde{\bf Dih}_{S}}
(\ul s,g\cdot \ul s)$ for $(t,\ul s)\in \on{Ob}\widetilde{\bf PaDih}_{S}$. 
We then set ${\bf PaDih}_{S}:= \widetilde{\bf PaDih}_{S}/D$. 

We have $\on{Ob}{\bf PaDih}_{S} = \{($a planar 3-valent tree, a map $\{$leaves$\}\to 
S)\}/($mirror symmetry$) = \sqcup_{n\geq 0}(PlT_{n}\times S^{n})/D_{n}$. 
We define a functor ${\bf PaB}_{S}^{hbal}\to {\bf PaDih}_{S}$ by the condition 
that: (a) at the level of objects, it is given by the canonical map 
$T_{n}\times S^{n}\to (PlT_{n}\times S^{n})/D_{n}$, (b) the diagram 
$\begin{matrix} {\bf PaB}_{S}^{hbal}& \to & {\bf PaDih}_{S}\\
\downarrow  & & \downarrow \\
{\bf B}_{S}^{hbal}  & \to &{\bf Dih}_{S}\end{matrix}$
commutes. One proves as above that this is a half-balanced contraction. 
Using the arguments of the proofs of Propositions \ref{prop:univ:tr}, \ref{bmc}
and \ref{prop:hrbmc}, one proves: 

\begin{proposition} \label{prop:uniqueness}
Let $\cC$ be a half-balanced braided monoidal category, equipped with a map
$S\stackrel{f}{\to} \on{Ob}\cC$ such that $f(s)^{*} = f(s)$ for any $s\in S$, 
and a balanced contraction $\cC\to\cO$. Then there exists a functor ${\bf PaDih}_{S}\to
\cO$, such that the diagram 
$\begin{matrix} {\bf PaB}_{S} & \to & {\bf PaDih}_{S} \\
 \downarrow &&\downarrow \\
 {\bf B}_{S}^{hbal} &\to&{\bf Dih}_{S}\end{matrix}$ commutes. 
\end{proposition}

We then have natural diagrams 
$$
\begin{matrix}
{\bf B}_{S} &\to &{\bf B}_{S}^{bal} &\to &{\bf B}_{S}^{hbal} \\
 & & \downarrow & & \downarrow \\
 & & {\bf Cyc}_{S} & \to & {\bf Dih}_{S}
 \end{matrix} \quad \on{and} \quad 
\begin{matrix}
{\bf PaB}_{S} &\to &{\bf PaB}_{S}^{bal} &\to &{\bf PaB}_{S}^{hbal} \\
 & & \downarrow & & \downarrow \\
 & & {\bf PaCyc}_{S} & \to & {\bf PaDih}_{S}
 \end{matrix} 
 $$
These diagrams fit in a bigger diagram, with a collection of functors from the left to the 
right-hand side diagram. 

\section{Completions} \label{sec:compl}

Let $G\to S_{n}$ be a group morphism. One can define the relative pro-$l$
and relative prounipotent completions $G_{l}$ and $G(-)$ of $G\to S_{n}$. They 
fit in exact sequences $1\to U_{l}\to G_{l}\to S_{n}\to 1$ and 
$1\to U(-)\to G(-)\to S_{n}\to 1$, where $U_{l}$ and $U(-)$ are 
pro-$l$ and $\QQ$-prounipotent. We have a morphism $G_{l}\to G(\QQ_{l})$
(\cite{HM}, Lemma A.7), fitting in a sequence of morphisms
$G\to \widehat G\to G_{l}\to G(\QQ_{l})$, where $\widehat G$ is the 
profinite completion of $G$. 
Applying this to $B_{n}$ are any of this quotients $B_n/Z_n$, $\Gamma_{0,n}$
considered above, we obtain for each of the categories $\cC = {\bf (Pa)B}_{S}^{(h)(bal)}$,
${\bf (Pa)Cyc}_{S}$, ${\bf (Pa)Dih}_{S}$, completed categories 
$\widehat \cC$, $\cC_{l}$, $\cC(-)$, and functors $\cC\to\widehat\cC\to
\cC_{l}\to \cC(\QQ_{l})$. 

Let us say that a pro-$l$ (resp., prounipotent) b.m.c.~is a b.m.c.~$\cC$, equipped with 
an assignment $\on{Ob}\cC\ni X\mapsto \cU_{X}\triangleleft \on{Aut}_{\cC}(X)$, 
such that $\cU_{X}$ is pro-$l$ (resp., prounipotent) for any $X$, and
for any $X,Y\in\on{Ob}\cC$ and $f\in\on{Iso}_{\cC}(X,Y)$, 
$f\cU_{X}f^{-1} = \cU_{Y}$ and $\on{im}(P_{n}\to \on{Aut}_{\cC}
(X_{1}\otimes\cdots\otimes X_{n})) \subset 
\cU_{X_{1}\otimes\cdots\otimes X_{n}}$ (here $P_{n} = 
\on{Ker}(B_{n}\to S_{n})$ is the pure braid group 
with $n$ strands). Similarly, $\cC$ is called profinite if $\on{Aut}_{\cC}(X)$
is profinite for any $X\in\on{Ob}\cC$. 

Then the completions $\widehat{\bf (Pa)B}_{S}$, ${\bf (Pa)B}_{S,l}$
and ${\bf (Pa)B}_{S}(-)$ are profinite, pro-$l$ and prounipotent (strict)
braided monoidal categories and are universal for such braided monoidal categories
$\cC$, equipped with a map $S\to\on{Ob}\cC$. 

\section{Actions of the Grothendieck--Teichm\"uller group}

\subsection{Grothendieck-Teichm\"uller semigroups}
Recall that the Grothendieck--Teichm\"uller semigroup is defined (\cite{Dr})
as 
\begin{eqnarray*}
&& \ul{\on{GT}} = \{(\lambda,f)\in (1+2\ZZ)\times F_{2} | 
f(Y,X)=f(X,Y)^{-1}, \\
 && f(X_{3},X_{1})X_{3}^{m}f(X_{2},X_{3})X_{2}^{m}
f(X_{1},X_{2})X_{1}^{m}=1, \quad \partial_{3}(f)\partial_{1}(f) = \partial_{0}(f)
\partial_{2}(f)\partial_{4}(f)\}, 
\end{eqnarray*}
where $F_{2}$ is the free group with two generators $X,Y$, 
$\partial_{0},\ldots,\partial_{4} : F_{2}\to P_{4}$ are simplicial 
morphisms, $X_{1}X_{2}X_{3}=1$, $m=(\lambda-1)/2$. It is a semigroup 
with $(\lambda,f)(\lambda',f')=(\lambda'',f'')$, where $\lambda'' = 
\lambda\lambda'$ and $f'' = \theta_{(\lambda',f')}(f)f'$, where 
$\theta_{(\lambda',f')}\in\on{End}(F_{2})$ is given by 
$(X,Y)\mapsto (f'X^{\lambda'}f^{\prime-1}, Y^{\lambda'})$. 
Then $\ul{\on{GT}}\to \on{End}(F_{2})^{op}$, $(\lambda,f)
\mapsto \theta_{(\lambda,f)}$ is a semigroup morphism. 
The profinite, pro-$l$ and prounipotent analogues 
$\widehat{\ul{\on{GT}}}$, $\ul{\on{GT}}_{l}$ and 
$\ul{\on{GT}}(-)$ of $\ul{\on{GT}}$ are defined by replacing 
$(\ZZ,F_{2})$ by $(\widehat\ZZ,\widehat F_{2})$, 
$(\ZZ_{l},(F_{2})_{l})$, and $\kk\mapsto (\kk,F_{2}(\kk))$
where $\kk$ is a $\QQ$-ring. We then have morphisms of semigroups
$\ul{\on{GT}}\to \widehat{\ul{\on{GT}}}\to 
\ul{\on{GT}}_{l}\to\ul{\on{GT}}(\QQ_{l})$; the associated 
groups are denoted $\on{GT},\widehat{\on{GT}},\on{GT}_{l},
\on{GT}(-)$.

\subsection{Action on (half-)braided monoidal categories}
The semigroup $\ul{\on{GT}}$ acts on $\{$braided monoidal categories$\}$
as follows: $(\lambda,f)*(\cC,\otimes,\beta_{XY},a_{XYZ}) = 
(\cC,\otimes,\beta'_{XY},a'_{XYZ})$, where $\beta'_{XY} = 
\beta_{XY}(\beta_{YX}\beta_{XY})^{m}$ and 
$$a'_{XYZ} = a_{XYZ}f(\beta_{YX}\beta_{XY}\otimes\on{id}_{Z},
a_{XYZ}^{-1}(\on{id}_{X}\otimes\beta_{ZY}\beta_{YZ})a_{XYZ}).$$ 
In the same way, $\widehat{\ul{\on{GT}}}$ acts on $\{$braided monoidal 
categories $\cC$, such that $\on{Aut}_{\cC}(X)$ is finite for any 
$X\in\on{Ob}\cC\}$, $\ul{\on{GT}}_{l}$ acts on $\{$pro-$l$ braided monoidal 
categories$\}$ and $\ul{\on{GT}}(\kk)$ acts on $\{\kk$-prounipotent 
braided monoidal categories$\}$. 

We have natural functors $\{$half-balanced braided monoidal categories$\}
\to\{$balanced braided monoidal categories$\}\to\{$braided monoidal categories$\}$. 

\begin{proposition}
The action of $\underline{\on{GT}}$ on $\{$braided monoidal categories$\}$
lifts to compatible actions on $\{$(half-)balanced braided monoidal categories$\}$. 
Similarly, the actions of $\widehat{\ul{\on{GT}}},\ldots,
\ul{\on{GT}}(\kk)$ lift to compatible actions on $\{$(half-)balanced finite braided 
monoidal categories$\}$, ...,  $\{$(half-)balanced $\kk$-prounipotent braided 
monoidal categories$\}$. 
\end{proposition}

{\em Proof.} This lift is given by $(\lambda,f)*(\cC,\otimes,\beta_{XY},a_{XYZ},
\theta_{X}):= (\cC,\otimes,\beta'_{XY},a'_{XYZ},\theta'_{X})$, 
where $\theta'_{X}:= \theta_{X}^{\lambda}$ and  
$(\lambda,f)*(\cC,\otimes,\beta_{XY},a_{XYZ},
a_{X}):= (\cC,\otimes,\beta'_{XY},a'_{XYZ},a'_{X})$, 
where $a'_{X}:= (a_{X^{*}}a_{X})^{m}a_{X}$, where $m=(\lambda-1)/2$.  
\hfill \qed\medskip 

\begin{proposition} \label{prop:comp:tr}
Let $\cC$ be a half-balanced category and let $\cC\stackrel{\langle-\rangle}{\to}
\cO$ be a half-balanced contraction. Then for any $(\lambda,f)\in\ul{\on{GT}}$,
the composite functor $(\lambda,f)*\cC\stackrel{\sim}{\to}\cC
\stackrel{\langle-\rangle}{\to}\cO$ is a half-balanced contraction on 
$(\lambda,f)*\cC$. Here $(\lambda,f)*\cC\stackrel{\sim}{\to}\cC$
is the identity functor (which is not tensor). Same statements with
$\cC$ finite, ..., $\kk$-unipotent and $\ul{\on{GT}}$ replaced by 
$\widehat{\ul{\on{GT}}}$, ..., $\ul{\on{GT}}(\kk)$. 
\end{proposition}

{\em Proof.} Assume that $(\cC,\beta_{XY},a_{X})$ is half-balanced; we set 
$\theta_{X}:= a_{X^{*}}a_{X}$. Then $(\cC,\beta_{XY},\theta_{X})$
is balanced and $\theta'_{X} = \theta_{X}^{\lambda}$.  Then $(\theta'_{Y}
\otimes \on{id}_{X})\beta'_{XY} = (\theta_{Y}
\otimes \on{id}_{X})\beta_{XY} (\theta_{X}^{-1}\otimes
\theta_{Y})^{m}\theta_{X\otimes Y}^{m}$. The identities 
$\langle\theta_{X}\rangle = \on{id}_{\langle X \rangle}$, 
$\langle \theta_{X}^{-1}\otimes \theta_{Y}\rangle = 
\on{id}_{\langle X,Y\rangle}$ (see Lemma \ref{lemma:hbal}) and 
$\langle (\theta_{Y}
\otimes \on{id}_{X})\beta_{XY} \rangle = \on{id}_{\langle X,Y\rangle}$
(as $\langle-\rangle$ is a half-balanced contraction) imply that 
$\langle (\theta'_{Y}
\otimes \on{id}_{X})\beta'_{XY} \rangle = \on{id}_{\langle X,Y\rangle}$, 
so $\langle-\rangle$ is a balanced contraction for $(\lambda,f)*\cC$. Moreover, 
$a'_{X} = a_{X}(a_{X^{*}}a_{X})^{m} = a_{X}\theta_{X}^{m}$, so 
$\langle\theta_{X}\rangle = \on{id}_{\langle X\rangle}$ implies 
$\langle a'_{X}\rangle = \langle a_{X}\rangle = \on{id}_{\langle X\rangle}$. 
\hfill \qed\medskip 

\subsection{Action on ${\bf PaDih}_{S}$}

For $(\lambda,f)\in\ul{\on{GT}}$, let $i_{(\lambda,f)}$ be the endomorphism of 
${\bf PaB}_{S}^{(h)bal}$ defined as the composite functor 
${\bf PaB}_{S}^{{(h)bal}}\stackrel{\alpha_{(\lambda,f)}}{\to}
(\lambda,f)*{\bf PaB}_{S}^{{(h)bal}}\stackrel{\sim}{\to}
{\bf PaB}_{S}^{(h)bal}$, where the first functor is the unique (half-)balanced 
monoidal functor which is the identity on objects, and the second functor is the 
identity functor (which is not monoidal). As in \cite{E}, Proposition 80, 
one shows that $(\lambda,f)\mapsto i_{(\lambda,f)}$ is a morphism 
$\ul{\on{GT}}\to \on{End}({\bf PaB}_{S}^{(h)bal})^{op}$. One 
similarly defines morphisms $\widehat{\ul{\on{GT}}}\to \on{End}
(\widehat{\bf PaB}_{S}^{(h)bal})^{op}$, ..., 
$\ul{\on{GT}}(\kk)\to \on{End}({\bf PaB}_{S,\kk}^{(h)bal})^{op}$. 

For $(\lambda,f)\in\ul{\on{GT}}$, we define an endofunctor 
$j_{(\lambda,f)}$ of ${\bf PaDih}_{S}$ as follows: according to 
Proposition \ref{prop:comp:tr}, the composite functor 
$(\lambda,f)*{\bf PaB}_{S}^{hbal}\stackrel{\sim}{\to}
{\bf PaB}_{S}^{hbal}\stackrel{\langle-\rangle}{\to}{\bf PaDih}_{S}$
 is a half-balanced contraction. By universality of the contraction ${\bf PaB}_{S}^{hbal}
\stackrel{\langle-\rangle}{\to}  {\bf PaDih}_{S}$, there exists a unique 
endofunctor $j_{(\lambda,f)}$  of ${\bf PaDih}_{S}$, such that 
the following diagram commutes
$$
\xymatrix{
{\bf PaB}_S^{hbal} \ar[r]_{\alpha_{(\lambda,f)}}\ar[d]_{\langle-\rangle}
\ar@/^{1pc}/[rr]^{i_{(\lambda,f)}} & (\lambda,f)*{\bf PaB}_S^{hbal} 
\ar[r]^{\sim}\ar[d] & {\bf PaB}_S^{hbal}\ar[dl]^{\langle-\rangle} \\
{\bf PaDih}_S \ar[r]_{j_{(\lambda,f)}} & {\bf PaDih}_S 
}
$$

\begin{proposition}
The map $(\lambda,f)\mapsto j_{(\lambda,f)}$ defines a morphism 
$\ul{\on{GT}}\to\on{End}({\bf PaDih}_{S})^{op}$; one similarly defines 
morphisms 
$\widehat{\ul{\on{GT}}}\to\on{End}(\widehat{\bf PaDih}_{S})^{op}$, 
etc.  
\end{proposition} 

{\em Proof.} We have a commutative diagram 
\begin{equation} \label{diagram:1}
\xymatrix{
{\bf PaB}_S^{hbal} \ar[r]^{\alpha_{(\lambda',f')}} \ar[dd]_{\langle-\rangle} & 
(\lambda',f')*{\bf PaB}_S^{hbal} \ar[d]^{\sim} \\
 & {\bf PaB}_S^{hbal} \ar[d]^{\langle-\rangle} \\
 {\bf PaDih}_S \ar[r]_{j_{(\lambda',f')}} & {\bf PaDih}_S} 
\end{equation}
which gives rise to 
$$
\xymatrix{
(\lambda,f)*{\bf PaB}_S^{hbal} \ar[r]^{(\lambda,f)*\alpha_{(\lambda',f')}} 
\ar[dd]_{\langle-\rangle} & 
(\lambda,f)(\lambda',f')*{\bf PaB}_S^{hbal} \ar[d]^{\sim} \\
 & {\bf PaB}_S^{hbal} \ar[d]^{\langle-\rangle} \\
 {\bf PaDih}_S \ar[r]_{j_{(\lambda',f')}} & {\bf PaDih}_S
} 
$$
Composing it with the analogue of (\ref{diagram:1}) with $(\lambda',f')$ 
replaced by $(\lambda,f)$, we get a commutative diagram 
$$
\xymatrix{
{\bf PaB}_S^{hbal} \ar[rrr]^{((\lambda,f)*\alpha_{(\lambda',f')}) 
\circ \alpha_{(\lambda,f)}} 
\ar[dd]_{\langle-\rangle} & &&
(\lambda,f)(\lambda',f')*{\bf PaB}_S^{hbal} \ar[d]^{\sim} \\
 & &&{\bf PaB}_S^{hbal} \ar[d]^{\langle-\rangle} \\
 {\bf PaDih}_S \ar[rrr]_{j_{(\lambda',f')}\circ 
 j_{(\lambda,f)}} & &&{\bf PaDih}_S
} 
$$

On the other hand, both $((\lambda,f)*
\alpha_{(\lambda',f')})\alpha_{(\lambda,f)}$ and 
$\alpha_{(\lambda,f)(\lambda',f')}$ are tensor functors 
${\bf PaB}_{S}^{hbal}\to (\lambda,f)(\lambda',f')*{\bf PaB}_{S}^{hbal}$
of half-balanced braided monoidal categories, inducing the identity at the level 
of objects, and by the uniqueness of such functors, they coincide. 
The above diagram may therefore be rewritten as 
$$
\xymatrix{
{\bf PaB}_S^{hbal} \ar[rr]^{\alpha_{(\lambda,f)(\lambda',f')} } 
\ar[dd]_{\langle-\rangle}  &&
(\lambda,f)(\lambda',f')*{\bf PaB}_S^{hbal} \ar[d]^{\sim} \\
  &&{\bf PaB}_S^{hbal} \ar[d]^{\langle-\rangle} \\
 {\bf PaDih}_S \ar[rr]_{j_{(\lambda',f')}\circ 
 j_{(\lambda,f)}} &&{\bf PaDih}_S
} 
$$
which may be viewed as a functor between half-balanced categories 
with a contraction. 

On the other hand, another such a functor is 
$$
\xymatrix{
{\bf PaB}_S^{hbal} \ar[rr]^{\alpha_{(\lambda,f)(\lambda',f')} } 
\ar[dd]_{\langle-\rangle}  &&
(\lambda,f)(\lambda',f')*{\bf PaB}_S^{hbal} \ar[d]^{\sim} \\
  &&{\bf PaB}_S^{hbal} \ar[d]^{\langle-\rangle} \\
 {\bf PaDih}_S \ar[rr]_{j_{(\lambda,f)(\lambda',f')}} &&{\bf PaDih}_S
} 
$$
By the universality of the contraction ${\bf PaB}_{S}^{hbal}
\stackrel{\langle-\rangle}{\to}{\bf PaDih}_S$, we then have 
$j_{(\lambda,f)(\lambda',f')} = j_{(\lambda',f')}j_{(\lambda,f)}$. 
\hfill \qed\medskip 

\subsection{Action on Teichm\"uller groupoids and proof of Theorem \ref{thm:1}}

$T_{0,S}$ may be viewed as the full subcategory of ${\bf PaDih}_{S}$
whose objects are the classes modulo $D$ of $PlT_{|S|}\times 
\on{Bij}(|S|,S)$. The action of $\ul{\on{GT}}$ then restricts to 
$T_{0,S}$, and similarly in the completed cases. In the profinite case, one
checks that that resulting action coincides with that defined in 
 in \cite{Sch}. This proves Theorem \ref{thm:1}

\subsection{Proof of Theorem \ref{thm:2}}

We define $T_{0,n}(\kk)$ by $\on{Ob}T_{0,n}(\kk) = \on{Ob}T_{0,n}$ 
and for $b,c\in\on{Ob}T_{0,n}$, $\on{Hom}_{T_{0,n}(\kk)}(b,c) = 
\on{Aut}_{T_{0,n}}(b)(\kk)\times_{\on{Aut}_{T_{0,n}}(b)}
\on{Hom}_{T_{0,n}}(b,c)$, where for $G$ a finitely generated group, 
$G(\kk)$ is its prounipotent completion. 

If $\pi$ is a finitely generated group, we define the group scheme $\ul{\on{Aut}\pi}(-)$
by $\ul{\on{Aut}\pi}(\kk):= \on{Aut}((\on{Lie}\pi)^{\kk})$, where for 
$\on{Lie}\pi$ is the Lie algebra of the prounipotent completion of $\pi$, 
$\g^{\kk} = \on{lim}_{\leftarrow}(\g/\g_{n})\otimes\kk$, and
$\g_{0}=\g$, $\g_{n+1}=[\g,\g_{n}]$. We then have a morphism 
$\ul{\on{Aut}\pi}(\kk)\to \on{Aut}(\pi(\kk))$, $\theta\mapsto\theta_{*}$. 
$\on{Aut}(\pi_{l},\pi(\QQ_{l}))$
is then defined as $\{(\theta,\theta_{l})\in \ul{\on{Aut}\pi}(\QQ_{l})
\times \on{Aut}(\pi_{l}) | \theta_{*}i = i \theta_{l}\}$, where 
$i$ is the morphism $\pi_{l}\to\pi(\QQ_{l})$. 

If $G$ is a groupoid such that 
$\on{Iso}_{G}(b,c)\neq\emptyset$ for any $b,c\in\on{Ob}G$, then the choice
of $b\in\on{Ob}G$ gives rise to an isomorphism $\on{Aut}G\simeq \pi^{
\on{Ob}G - \{b\}}\rtimes \on{Aut}\pi$, where $\pi = \on{Aut}_{G}(b)$; we then 
define the group scheme $\ul{\on{Aut}G}(-)$ by $\ul{\on{Aut}G}(\kk):= 
\pi(\kk)^{\on{Ob}G - \{b\}}\rtimes\ul{\on{Aut}\pi}(\kk)$. We define 
as above $\on{Aut}(G_{l},G(\QQ_{l}))$ and the morphisms 
$\on{Aut}(G_{l})\leftarrow \on{Aut}(G_{l},G(\QQ_{l})) \to 
\ul{\on{Aut}G}(\QQ_{l})$. 

We have morphisms $G_{\QQ}\to \on{GT}_{l}\to \on{GT}(\QQ_{l})$ and a 
functor ${\bf PaB}_{S,l}\to {\bf PaB}_{S}(\QQ_{l})$. Theorem \ref{thm:2}
follows from the fact that this functor is compatible with the actions of
$\on{GT}_{l},\on{GT}(\QQ_{l})$ on ${\bf PaB}_{S,l}$, 
${\bf PaB}_{S}(\QQ_{l})$.

\section{Graded aspects}

Let $\t_{n}$ be the graded Lie algebra with generators $t_{ij}$, 
$i\neq j\in [n]$ and relations $t_{ji}=t_{ij}$, $[t_{ij},t_{ik}+t_{jk}]=0$, 
$[t_{ij},t_{kl}]=0$ for $i,j,k,l$ distinct. Let $\p_{n}$ be the quotient of  
$\t_{n}$ by the relations $\sum_{j|j\neq i}t_{ij}=0$, for any $i\in [n]$. 
Equivalently, $\p_{n}$ is presented by generators $t_{ij}$ are relations 
$t_{ji}=t_{ij}$, $\sum_{j|j\neq i}t_{ij}=0$ for any $i$, and $[t_{ij},t_{kl}]=0$
for $i,j,k,l$ distinct. 

Let $\kk$ be a $\QQ$-ring, then the set $M(\kk)$ of Drinfeld associators 
defined over $\kk$ is the set of pairs $(\mu,\Phi)\in\kk\times
\on{exp}(\hat\f_{2}^{\kk})$, satisfying the  duality, hexagon and pentagon 
conditions\footnote{If $\g$ is a graded Lie algebra, then 
$\hat\g^{\kk}$ is the degree completion of $\g\otimes \kk$.} 
(see \cite{Dr}). The data of $t\in T_{n}$ and $(\mu,\Phi)\in 
M(\kk)$ gives rise to a morphism $B_{n}\stackrel{i_{t,\Phi}}{\to} 
\on{exp}(\hat\t_{n}^{\kk})\rtimes S_{n}$, which extends to an 
isomorphism $B_{n}(\kk)\stackrel{\sim}{\to} \on{exp}(\hat\t_{n}^{\kk})
\rtimes S_{n}$ (see e.g. \cite{AET}) if $\mu\in\kk^{\times}$. 

\begin{proposition}
There exists a unique morphism $\Gamma_{0,n}\to 
\on{exp}(\hat\p_{n}^{\kk})\rtimes S_{n}$, 
such that the diagram 
$\begin{matrix} B_{n} & \stackrel{i_{t,\Phi}}{\to} & 
\on{exp}(\hat\t_{n}^{\kk})\rtimes S_{n} \\
\downarrow & & \downarrow \\
\Gamma_{0,n} & \to & 
\on{exp}(\hat\p_{n}^{\kk})\rtimes S_{n}\end{matrix}$
commutes. It gives rise to an isomorphism $\Gamma_{0,n}(\kk)
\stackrel{\sim}{\to}\on{exp}(\hat\p_{n}^{\kk})\rtimes S_{n}$.  
\end{proposition}

{\em Proof.} One checks that $i_{t,\Phi}$ takes $z_{n}$ to 
$\on{exp}(\mu\sum_{1\leq i<j\leq n}t_{ij})$ and $\sigma_{i}\cdots 
\sigma_{n-1}^{2}\cdots \sigma_{0}^{2}\cdots \sigma_{i-1}$ to a 
conjugate of $\on{exp}(\mu\sum_{j|j\neq i}t_{ij})$. This implies 
the announced commutative diagram. Let $\Gamma_{0,[n]}:= \on{Ker}(
\Gamma_{0,n}\to S_{n})$ and $\Gamma_{0,[n]}(\kk)$ be its $\kk$-prounipotent
completion. The morphism $\Gamma_{0,n}\to \on{exp}(\hat\p_{n}^{\kk})\rtimes
S_{n}$ gives rise to a morphism $\Gamma_{0,[n]}(\kk)\to \on{exp}(\hat\p_{n}^{\kk})$;
let us show that this is an isomorphism. We have a morphism $\t_{n}\to \on{gr}
\on{Lie}P_{n}$, where $P_{n}:= \on{Ker}(B_{n}\to S_{n})$, given by 
$t_{ij}\mapsto$ class of $\on{log} (\sigma_{i}\cdots\sigma_{j-2})\sigma_{j-1}^{2}
 (\sigma_{i}\cdots\sigma_{j-2})^{-1}$. We then have a commutative diagram 
 $\begin{matrix}\t_{n} & \to & \on{gr}\on{Lie}P_{n} \\
 \downarrow & & \downarrow \\
 \p_{n} & \to & \on{gr}\on{Lie}\Gamma_{0,[n]}\end{matrix}$
 where the horizontal maps are surjective and the Lie algebras in the right 
 side are generated in degree 1. The Lie algebra morphism corresponding to the 
 group morphism $\Gamma_{0,[n]}(\kk)\to \on{exp}(\hat\p_{n}^{\kk})$
 is a Lie algebra morphism $\on{Lie}\Gamma_{0,[n]}(\kk)\to \hat\p_{n}^{\kk}$, 
 whose associated graded morphism is a graded Lie algebra morphism 
 $\on{gr}\on{Lie}\Gamma_{0,[n]}(\kk)\to \p_{n}^{\kk}$. The composite 
 map $\p_{n}^{\kk}\to \on{gr}\on{Lie}\Gamma_{0,[n]}\otimes \kk
 \to \p_{n}^{\kk}$ is a graded isomorphism, as it can be checked on the 
 degree 1 part of $\p_{n}$. It follows that the morphism $\p_{n}\to \on{gr}
 \on{Lie}\Gamma_{0,[n]}$ is injective as well, therefore it is an isomorphism 
of Lie algebras. So $\Gamma_{0,[n]}(\kk)\to \on{exp}(\hat\p_{n}^{\kk})$
is an isomorphism. \hfill \qed\medskip
 
We define a category ${\bf PaDih}_{S}^{gr}$ similarly to ${\bf PaDih}_{S}$, 
i.e., as the quotient by $D$ of an intermediate category $\widetilde{\bf PaDih}_{S}^{gr}$
obtained from $\widetilde{\bf PaDih}_{S}$ by replacing $\Gamma_{0,n}$ by 
$\on{exp}(\hat\p_{n}^{\kk})\rtimes S_{n}$, 
and the morphism $D\to D_{n}\to \Gamma_{0,n}$ by $D\to D_{n}\to S_{n}\to
\on{exp}(\hat\p_{n}^{\kk})\rtimes S_{n}$. 

If $(\mu,\Phi)\in M(\kk)$, recall that a braided monoidal category 
${\bf PaCD}^{\Phi}_{S}$ may be defined as follows: 
$\on{Ob}{\bf PaCD}^{\Phi}_{S} = \on{Ob}{\bf PaB}_{S}$; 
$\on{Hom}_{{\bf PaCD}^{\Phi}_{S}}((\ul s,t),(\ul s',t'))$
is empty if $|\ul s|\neq |\ul s'|$, and is equal to 
$\on{exp}(\hat\t_{n}^{\kk})\rtimes\{f\in S_{n}|\ul s'f = \ul s\}$; the 
composition is induced by the product in $\on{exp}(\hat\t_{n}^{\kk})\rtimes
S_{n}$; and the tensor product is obtained by restriction from the group morphism 
$(\on{exp}(\hat\t_{n}^{\kk})\rtimes S_{n})\times 
(\on{exp}(\hat\t_{n'}^{\kk})\rtimes S_{n'})\to 
\on{exp}(\hat\t_{n+n'}^{\kk})\rtimes S_{n+n'}$, induced by the Lie 
algebra morphism $\hat\t_{n}^{\kk}\times\hat\t_{n'}^{\kk}\to 
\hat\t_{n+n'}^{\kk}$, $(t_{ij},0)\mapsto t_{ij}$, $(0,t_{ij})\mapsto t_{n+i,n+j}$, 
and the group morphism $S_{n}\times S_{n'}\to S_{n+n'}$, $(\sigma,\sigma')\mapsto 
\sigma*\sigma'$, such that $(\sigma*\sigma')(i) = \sigma(i)$ for $i\in [n]$, and 
$(\sigma*\sigma')(n+i) = n+\sigma'(i)$ for $i\in [n']$. The braiding constraint
is defined by $\beta_{XY} = (e^{\mu t_{12}/2})^{[n],n+[n']}s_{n,n'}$ and the
associativity constraint is defined by $a_{XYZ} = (\Phi(t_{12},t_{23}))^{[n],n+[n'],
n+n'+[n'']}$ for $|X|=n$, $|Y|=n'$, $|Z| = n''$, $s_{n,n'}\in S_{n+n'}$ is defined by 
$s_{n,n'}(i)=n'+i$ for $i\in [n]$ and $s_{n,n'}(n+i)=i$ for $i\in [n']$, and 
for $I_{1},\ldots,I_{n}\subset [m]$ disjoint subsets, the morphism 
$\t_{n}\to \t_{m}$, $x\mapsto x^{I_{1},\dots,I_{n}}$ is defined by 
$t_{ij}\mapsto \sum_{\alpha\in I_{i},\beta\in I_{j}}t_{\alpha\beta}$. 
Then ${\bf PaCD}_{S}^{\Phi}$ is a braided monoidal category; it follows that there
is a unique monoidal functor $j_\Phi:{\bf PaB}_{S}\to {\bf PaCD}_{S}^{\Phi}$, 
which induces the identity on objects. 

We then define a functor ${\bf PaCD}_{S}^{\Phi}\to{\bf PaDih}_{S}^{gr}$
as the composite functor ${\bf PaCD}_{S}^{\Phi}\to
\widetilde{\bf PaDih}_{S}^{gr}\to{\bf PaDih}_{S}^{gr}$, where the first functor 
is induced by the projection morphisms $\t_{n}\to\p_{n}$ and the second functor 
is the quotient functor $\widetilde{\bf PaDih}_{S}^{gr}\to
\widetilde{\bf PaDih}_{S}^{gr}/D \simeq {\bf PaDih}_{S}^{gr}$. 

\begin{proposition} \label{prop:hbal}
The functor ${\bf PaCD}_{S}^{\Phi}\to{\bf PaDih}_{S}^{gr}$ is a half-balanced 
contraction. 
\end{proposition}

{\em Proof.} We first show: 

\begin{lemma}
Let $X\in\on{Ob}{\bf PaB}_{\{\bullet\}}$ be of degree $n$, then 
$$\on{im}(h_{n}\in{\bf PaB}_{\{\bullet\}}(X,X^{*})\to 
{\bf PaCD}_{\{\bullet\}}^{\Phi}(X,X^{*})) = \on{exp}({\mu\over 2}
\sum_{1\leq i<j\leq n}t_{ij})\bigl(\begin{smallmatrix} 1 & 2 & \ldots & n \\
n & n-1 & \ldots & 1\end{smallmatrix}\bigl)\in \on{exp}(\hat\t_{n}^{\kk})
\rtimes S_{n}.$$ 
\end{lemma}

{\em Proof.} It suffices to prove this for a particular $X_{0}\in
\on{Ob}{\bf PaB}_{\{\bullet\}}$ of degree $n$, say $X_{0} = \bullet
(\bullet(\cdots(\bullet\bullet)))$. Indeed, if we denote by 
$h_{n}^{X}\in {\bf PaB}_{\{\bullet\}}(X,X^{*})$ the element corresponding to 
$h_{n}$ and if we have $\on{im}(h_{n}^{X_{0}}) =  \on{exp}({\mu\over 2}
\sum_{1\leq i<j\leq n}t_{ij})\bigl(\begin{smallmatrix} 1 & 2 & \ldots & n \\
n & n-1 & \ldots & 1\end{smallmatrix}\bigl)$, then if $X$ is another object with the
same degree, then $\on{im}(h_{n}^{X}) = \Phi_{X_{0}^{*},X^{*}}
\on{im}(h_{n}^{X_{0}})\Phi_{X,X_{0}}$, where $\Phi_{X,Y} = 
\on{im}(1\in {\bf PaB}_{\{\bullet\}}(X,Y)\to {\bf PaCD}^{\Phi}_{
\{\bullet\}}(X,Y) = \on{exp}(\hat\t_{n}^{\kk})\rtimes S_{n})$. 
As $\Phi_{X_{0}^{*},X^{*}} = \bigl(\begin{smallmatrix} 1 & 2 & \ldots & n \\
n & n-1 & \ldots & 1\end{smallmatrix}\bigl) \Phi_{X_{0},X}
\bigl(\begin{smallmatrix} 1 & 2 & \ldots & n \\
n & n-1 & \ldots & 1\end{smallmatrix}\bigl)$, $\sum_{1\leq i<j\leq n}
t_{ij}\in\t_{n}$ is central and $\Phi_{X_{0}^{*},X^{*}} 
= \Phi_{X^{*},X_{0}^{*}}^{-1}$, $\on{im}(h_{n}^{X}) =  \on{exp}({\mu\over 2}
\sum_{1\leq i<j\leq n}t_{ij})\bigl(\begin{smallmatrix} 1 & 2 & \ldots & n \\
n & n-1 & \ldots & 1\end{smallmatrix}\bigl)$. 

We now prove that statement for $X_{0} = \bullet
(\bullet(\cdots(\bullet\bullet)))$ of degree $n$, which we redenote $X_{n}$. 
The proof is by induction on $n$. The statement is clear for $n=1,2$. Assume it at order 
$n-1$. Then $h_{n}^{X_{n}}\in{\bf PaB}_{\{\bullet\}}(X_{n},X_{n}^{*})$ may
be viewed as the composite morphism $X_{n} = \bullet \otimes X_{n-1}
\stackrel{\on{id}_{\bullet}\otimes h_{n-1}}{\to} \bullet\otimes
X_{n-1}^{*}\stackrel{\beta_{\bullet,X_{n-1}^{*}}}{\to} X_{n-1}^{*}
\otimes\bullet = X_{n}^{*}$, whose image in ${\bf PaCD}_{\{\bullet\}}^{\Phi}$
is $s \on{exp}({\mu\over 2}\sum_{i=2}^{n}t_{1i})\on{exp}({\mu\over 2}
\sum_{2\leq i<j\leq n}t_{ij})\bigl(\begin{smallmatrix} 1 & 2 & \ldots & n \\
1 & n & \ldots & 2\end{smallmatrix}\bigl)$, where $s = 
\bigl(\begin{smallmatrix} 1 & 2 & \ldots & n \\
n & 1 & \ldots & n-1\end{smallmatrix}\bigl)$. So this image is 
$  \on{exp}({\mu\over 2}
\sum_{1\leq i<j\leq n}t_{ij})\bigl(\begin{smallmatrix} 1 & 2 & \ldots & n \\
n & n-1 & \ldots & 1\end{smallmatrix}\bigl)$. \hfill \qed\medskip 

We then show: 

\begin{lemma}
Let $X,Y\in\on{Ob}{\bf PaB}_{\{\bullet\}}$ be of degrees $n,m$,
then 
\begin{eqnarray*}
&& \on{im}((\theta_{Y}\otimes\on{id}_{X})\beta_{XY}\in
{\bf PaB}_{\{\bullet\}}(X\otimes Y,Y\otimes X)\to
{\bf PaCD}_{\{\bullet\}}^{\Phi}(X\otimes Y,Y\otimes X))
\\ && = \bigl(\begin{smallmatrix} 1 & \cdots & n & n+1 & \cdots & n+m \\
m+1 & \cdots & m+n & 1 & \cdots & m\end{smallmatrix}\bigl)
\on{exp}({\mu\over 2}\sum_{j\in n+[m]}\sum_{\alpha\in[n+m] - \{j\}}
t_{\alpha j}).
\end{eqnarray*}  
\end{lemma}

{\em Proof.} The image of $\beta_{XY}$ is
$\bigl(\begin{smallmatrix} 1 & \cdots & n & n+1 & \cdots & n+m \\
m+1 & \cdots & m+n & 1 & \cdots & m\end{smallmatrix}\bigl)
\on{exp}({\mu\over 2}\sum_{i\in [n],j\in n+[m]}t_{ij})$, while the image of 
$\theta_{Y}\otimes \on{id}_{X}$ is $\on{exp}(\mu\sum_{j<j'\in [m]}
t_{jj'})$. \hfill \qed\medskip 

{\em End of proof of Proposition \ref{prop:hbal}.} If $X\in
\on{Ob}{\bf PaCD}_{S}^{\Phi}$ has degree $n$, then the image of 
$a_{X}\in {\bf PaCD}_{S}^{\Phi}(X,X^{*})$ in 
$\widetilde{\bf PaDih}_{S}^{gr}(X,X^{*}) = \on{exp}(\hat\p_{n}^{\kk})
\rtimes S_{n}$ is $\bigl(\begin{smallmatrix} 1 & 2 & \cdots & n \\
n & n-1 & \cdots & 1\end{smallmatrix}\bigl)$ as $\on{im}(\sum_{i<j\in [n]}
t_{ij}\in \t_{n}\to\p_{n})=0$. Now 
$\bigl(\begin{smallmatrix} 1 & 2 & \cdots & n \\
n & n-1 & \cdots & 1\end{smallmatrix}\bigl) = i_{X}^{s}$, therefore 
after taking the quotient by $D$, $\langle a_{X}\rangle = \on{id}_{\langle X\rangle}$
in $\on{End}_{{\bf PaDih}_{S}^{gr}}(\langle X\rangle)$. 

Similarly, if $X,Y\in\on{Ob}({\bf PaCD}_{S}^{\Phi})$ have degrees $n,m$, 
then the image of $(\theta_{Y}\otimes\on{id}_{X})\beta_{XY}\in
{\bf PaCD}_{S}^{\Phi}(X\otimes Y,Y\otimes X)$ in 
$\widetilde{\bf PaDih}_{S}^{gr}(X\otimes Y,Y\otimes X) = 
\on{exp}(\hat\p_{n+m}^{\kk})\rtimes S_{n+m}$ is $c^{m}$, 
where $c = \bigl(\begin{smallmatrix}1 & 2 & \cdots & n+m \\
2 & 3 & \cdots & 1\end{smallmatrix}\bigl)$ as $\on{im}(\sum_{\alpha\in
[n+m]-\{j\}}t_{\alpha j}\in\t_{n+m}\to\p_{n+m})=0$ for any $j$. 
It follows that this image coincides with $i_{X\otimes Y}^{r^{m}}$, whose
image in $\on{Aut}_{{\bf PaDih}_{S}^{gr}}(\langle X\otimes Y\rangle)$
is $\on{id}_{\langle X \otimes Y\rangle}$. It follows that 
$\langle(\theta_{Y}\otimes\on{id}_{X})\beta_{XY}\rangle
= \on{id}_{\langle X\otimes Y\rangle}\in\on{Aut}_{{\bf PaDih}_{S}^{gr}}
(\langle X\otimes Y\rangle)$. All this implies that ${\bf PaCD}_{S}^{\Phi}
\to{\bf PaDih}_{S}^{gr}$ satisfies the half-balanced contraction conditions. \hfill\qed\medskip

Proposition \ref{prop:hbal} immediately implies: 

\begin{corollary}
There exists a unique functor ${\bf PaDih}_S \stackrel{k_\Phi}{\to}
{\bf PaDih}_S^{gr}$, with is the identity on objects and such that 
the diagram $\xymatrix{{\bf PaB}_S \ar[r]^{j_\Phi}
\ar[d]^{\langle-\rangle}& {\bf PaCD}_S^\Phi\ar[d]^{\langle-\rangle}
\\ {\bf PaDih}_S \ar[r]^{k_\Phi} & {\bf PaDih}_S^{gr}}$
commutes. 
\end{corollary}

Recall that the graded Grothendieck-Teichm\"uller group $\on{GRT}(\kk)$ is 
defined as $\on{GRT}(\kk) = \on{GRT}_{1}(\kk)\rtimes\kk^{\times}$, 
where $\on{GRT}_1(\kk)$ is the set of all $g\in\on{exp}(\hat\f_2^\kk)
\subset\on{exp}(\hat\t_3^\kk)$ ($\f_2\subset\t_3$ being the Lie subalgebra
generated by $t_{12},t_{23}$), such that 
$$
g^{3,2,1}=g^{-1},\quad t_{12}+\on{Ad}(g^{1,2,3})^{-1}(t_{23})
+\on{Ad}(g^{2,1,3})^{-1}(t_{13})=t_{12}+t_{23}+t_{13}, 
$$
$$
g^{2,3,4}g^{1,23,4}g^{1,2,3}=g^{1,2,34}g^{12,3,4}, 
$$
equipped with the group law $(g_1*g_2)(A,B) := g_1(\on{Ad}(g_2(A,B))(A),B)
g_2(A,B)$, on which $\kk^\times$ acts by $(c\cdot g)(A,B)
:=g(c^{-1}A,c^{-1}B)$.

We now construct an action of this group on ${\bf PaDih}_{S}^{gr}$. 
For this, we recall from \cite{E} the notion of infinitesimally braided monoidal
category (i.b.m.c.). 

\begin{definition}
An i.b.m.c.~is a braided monoidal category $(\cC,\otimes,c_{XY},a_{XYZ})$, 
which is 

(1) symmetric, i.e., such that $c_{YX}c_{XY} = \on{id}_{X\otimes Y}$
for any $X,Y\in\on{Ob}\cC$, 

(2) prounipotent (see Section \ref{sec:compl}), i.e., equipped with an assignment 
$\on{Ob}\cC\ni X\mapsto \cU_{X}\triangleleft\on{Aut}_{\cC}(X)$, 
such that $f\cU_{X}f^{-1} = \cU_{Y}$ for $f\in\on{Iso}_{\cC}(X,Y)$, 

(3) equipped with a functorial assignment $(\on{Ob}\cC)^{2}\ni (X,Y)
\mapsto t_{XY}\in \on{Lie}\cU_{X\otimes Y}$, such that $t_{YX} = 
c_{YX}t_{YX}c_{XY}$ and 
$$
t_{X\otimes Y,Z} = a_{XYZ}(\on{id}_{X}\otimes t_{YZ})a_{XYZ}^{-1}
+(c_{YX}\otimes\on{id}_{Z})a_{YXZ}(\on{id}_{Y}\otimes t_{XZ})
a_{YXZ}^{-1} (c_{YX}\otimes\on{id}_{Z})^{-1}.$$
\end{definition}

According to \cite{Dr}, $\on{GRT}(\kk)$ acts on $\{$i.m.b.~categories$\}$ from the 
right as follows: $g\in \on{GRT}_{1}(\kk)\subset \on{exp}(\hat\f_{2}^{\kk})$
acts by $(\cC,\otimes,c_{XY},a_{XYZ},t_{XY})\cdot g := 
(\cC,\otimes,c_{XY},a'_{XYZ},t_{XY})$, where $a'_{XYZ}:= 
g(t_{XY}\otimes\on{id}_{Z},a_{XYZ}(\on{id}_{X}\otimes t_{YZ})
a_{XYZ}^{-1})a_{XYZ}$ and $c\in\kk^{\times}$ acts by 
$(\cC,\ldots)\cdot g := (\cC,\otimes,c_{XY},a_{XYZ},ct_{XY})$. 
Moreover, ${\bf PaCD}_{S}$, equipped with $c_{XY}:= s_{|X|,|Y|}$, 
$a_{XYZ}:= \on{id}_{|X|+|Y|+|Z|}$ and $t_{XY}:= t_{12}^{[|X|],|X|+[|Y|]}$
is universal among i.b.m.cs $\cC$, equipped with a map $S\to\on{Ob}\cC$. 
We derive from this, as in \cite{E}, Proposition 80, a morphism 
$\on{GRT}(\kk)\to \on{Aut}({\bf PaCD}_{S})$. 

We now introduce the notion of a balanced i.b.m.c.

\begin{definition}
A balanced structure on the i.b.m.c.~$\cC$ is a functorial assignment 
$\on{Ob}\cC\ni X\mapsto t_{X}\in\on{Lie}\cU_{X}$, such that 
for any $X,Y\in\on{Ob}\cC$, $t_{X\otimes Y} - t_{X}\otimes\on{id}_{Y}
- \on{id}_{X}\otimes t_{Y} = t_{XY}$. 
\end{definition}

\begin{definition}
A contraction on the small balanced i.b.m.c.~$\cC$ is a functor $\cC\stackrel{\langle-
\rangle}{\to}\cO$, such that for any $X,Y,Z\in\on{Ob}\cC$, 
$\langle X\otimes Y\rangle = \langle Y\otimes X\rangle(=:\langle
X,Y\rangle)$, $\langle (X\otimes Y)\otimes Z\rangle = \langle 
X\otimes (Y\otimes Z)\rangle (=:\langle X,Y,Z\rangle)$, 
$\langle c_{XY}\rangle = \on{id}_{\langle X,Y\rangle}$, 
$\langle a_{XYZ}\rangle = \on{id}_{\langle X,Y,Z\rangle}$, and 
$\langle t_{XY}+2\on{id}_{X}\otimes t_{Y}\rangle=0$. 
\end{definition}

\begin{remark}
We derive from the latter condition that $\langle t_{X}\rangle = 0$ for any 
$X\in\on{Ob}\cC$. Indeed, it gives by symmetrization $\langle t_{X\otimes Y}
\rangle =0$, and therefore $\langle t_{X}\rangle = 0$ by taking $Y={\bf 1}$. 
By antisymmetrization, this condition also implies $\langle t_{X}\otimes
\on{id}_{Y} - \on{id}_{X}\otimes t_{Y}\rangle = 0$. 
\end{remark}

We now construct a universal contraction on balanced i.b.m.~categories. 

\begin{proposition}
The i.b.m.c.~${\bf PaCD}_{S}$ is equipped with a balanced structure given by 
$t_{X} = \sum_{1\leq i<j\leq n}t_{ij}$ for $|X|=n$. Then the functor 
${\bf PaCD}_{S}\to{\bf PaDih}_{S}^{gr}$ is a contraction. 
\end{proposition}

{\em Proof.} For $|X|=n$, $|Y|=m$, $t_{XY}+2\on{id}_{X}\otimes t_{Y} = 
\sum_{j\in n+[m]}\sum_{\alpha\in[n+m]-\{j\}}t_{j\alpha}$, so 
$\langle t_{XY}+2\on{id}_{X}\otimes t_{Y}\rangle=0$ as
$\sum_{\alpha\in[n+m]-\{j\}}t_{j\alpha}=0$ in $\p_{n+m}$ for any 
$j\in n+[m]$. \hfill \qed\medskip 

\begin{proposition}
Let $\cC$ be a balanced i.b.m.c.~and let $\cC\stackrel{\langle-\rangle}{\to}
\cO$ be a contraction. Let $S\stackrel{f}{\to}\on{Ob}\cC$ be a map such that 
for any $s\in S$, $t_{f(s)}=0$. Then we have a commutative diagram 
$\begin{matrix} {\bf PaCD}_{S} &\to & \cC\\
\downarrow &&\downarrow \\
{\bf PaDih}_{S}^{gr}&\to&\cO\end{matrix}$
\end{proposition}

{\em Proof.} As $\cC$ is an i.b.m.c., there exists a unique 
functor ${\bf PaCD}_{S}\to\cC$ of i.b.m.~categories, extending $f$. 
As $t_{f(s)}=0$ for $s\in S$, it is compatible with the balanced structures. 
The construction of the commutative diagram is similar to the proof of 
Propositions \ref{bmc}, \ref{prop:hrbmc}. \hfill \qed\medskip 

\begin{proposition}
1) The action of $\on{GRT}(\kk)$ on $\{$i.b.m.~categories$\}$
lifts to $\{$balanced i.b.m.~categories$\}$ as follows: for
$(\cC,\otimes,c_{XY},a_{XYZ},t_{XY},t_{X})$ a balanced 
i.b.m.c., and $g\in \on{GRT}(\kk)$, 
$\cC\cdot g = (\cC,\ldots,t'_{X})$, where $t'_{X} = ct_{X}$
and $c = \on{im}(g\in \on{GRT}(\kk)\to\kk^{\times})$. 

2) If $\cC\stackrel{F}{\to}\cO$ is a contraction of the balanced i.b.m.c.~$\cC$, 
then $\cC\cdot g \stackrel{\sim}{\to}\cC\stackrel{F}{\to}\cO$
is a contraction of the balanced i.b.m.c.~$\cC$ (where $\cC\cdot g \stackrel{\sim}{\to}\cC$
is the identity of the underlying categories).  
\end{proposition}

The proof is immediate. 

We now construct an action of $\on{GRT}(\kk)$ on ${\bf PaCD}_{S}
\to{\bf PaDih}_{S}^{gr}$. A morphism $\on{GRT}(\kk)\to 
\on{Aut}({\bf PaCD}_{S})$, $g\mapsto a_{g}$ is defined by 
$a_{g} : {\bf PaCD}_{S}\to {\bf PaCD}_{S}*g \stackrel{\sim}{\to}
{\bf PaCD}_{S}$, where the first morphism is the unique functor of 
i.m.b.~categories, inducing the identity on objects, and the second 
morphism is the identification of the underlying categories.  

We define a morphism $\on{GRT}(\kk)\to\on{Aut}({\bf PaDih}_{S}^{gr})$, 
$g\mapsto j_{g}$ by the condition that the diagram 
$\begin{matrix}{\bf PaCD}_{S}&\to&{\bf PaCD}_{S}*g\\
\scriptstyle{\langle-\rangle}\downarrow&&\downarrow
\scriptstyle{\langle-\rangle}\\
{\bf PaDih}_{S}^{gr}&\stackrel{j_{g}}{\to}&{\bf PaDih}_{S}^{gr}\end{matrix}$
is a functor of balanced i.b.m.~categories with contractions. 
We then have a commutative diagram $\begin{matrix}{\bf PaCD}_{S}&
\stackrel{a_{g}}{\to}&{\bf PaCD}_{S}\\
\scriptstyle{\langle-\rangle}\downarrow&&
\downarrow\scriptstyle{\langle-\rangle}\\
{\bf PaDih}_{S}^{gr}&\stackrel{j_{g}}{\to}&{\bf PaDih}_{S}^{gr}
\end{matrix}$

Let now $(\mu,\Phi)\in M(\kk)$, where $\mu\in\kk^{\times}$, be an associator. 
It gives rise to an isomorphism $i_{\Phi}:\on{GT}(\kk)\to\on{GRT}(\kk)$, 
defined by the condition that $g*\Phi=\Phi*i_{\Phi}(g)$ for any $g\in\on{GT}(\kk)$. 
In the diagram 
$$
\xymatrix{
{\bf PaB}_S \ar[r]_{j_\Phi} \ar@/^{1pc}/[rr]^{\langle-\rangle} 
\ar[d]_{g}& {\bf PaCD}_S \ar[d]_{i_\Phi(g)} 
\ar@/^{1pc}/[rr]^{\langle-\rangle} & {\bf PaDih}_S
\ar[d]_{g} \ar[r]_{k_\Phi}& {\bf PaDih}_S^{gr}\ar[d]^{i_\Phi(g)}\\
{\bf PaB}_S\ar[r]^{j_\Phi} \ar@/_{1pc}/[rr]_{\langle-\rangle} &
{\bf PaCD}_S\ar@/_{1pc}/[rr]_{\langle-\rangle}  & 
{\bf PaDih}_S\ar[r]^{k_\Phi} & {\bf PaDih}_S^{gr}
}
$$
all the squares except perhaps the rightmost one commute. But this last 
square has to commute by the uniqueness of the morphism ${\bf PaDih}_{S}
\to\cO$ in Proposition \ref{prop:uniqueness} (the existence in this 
proposition implies uniqueness by abstract nonsense). 

All this implies that the isomorphism ${\bf PaDih}_{S}\stackrel{k_\Phi}{\to}
{\bf PaDih}_{S}^{gr}$ gives rise to a commutative diagram 
$\begin{matrix}
\on{GT}(\kk)&\to&\on{Aut}{\bf PaDih}_{S}\\
\scriptstyle{i_{\Phi}}{\downarrow}&&\downarrow\\
\on{GRT}(\kk)&\to&\on{Aut}{\bf PaDih}_{S}^{gr}
\end{matrix}$

The isomorphism ${\bf PaDih}_{S}\stackrel{k_\Phi}{\to}{\bf PaDih}_{S}^{gr}$
and the actions of $\on{G(R)T}(\kk)$ on these 
categories induce the identity at the level of objects. We then define
$T_{0,n}^{gr}$ to be the full subcategory of ${\bf PaDih}_{[n]}^{gr}$, 
whose set of objects is $(PlT_{n}\times \on{Bij}([n],[n]))/D_{n}$, and 
obtain this way an isomorphism $T_{0,n}(\kk)\to T_{0,n}^{gr}$
inducing a commutative diagram $\begin{matrix}\on{GT}(\kk)&\to&
\on{Aut}T_{0,n}(\kk)\\
\downarrow &&\downarrow\\
\on{GRT}(\kk)&\to&\on{Aut}T_{0,n}^{gr}(\kk)\end{matrix}$

This proves Theorem \ref{thm:3}. 

\begin{remark}
$T_{0,n}^{gr}$ could alternatively be defined as $\pi_{dih}^{*}\cC_{\Gamma,G,S}$,
where $\Gamma=D_{n}$, $G=\on{exp}(\hat\p_{n}^{\kk})\rtimes S_{n}$, and 
$S=[n]$ (see Section \ref{sec:1}). 
\end{remark}

\end{document}